\documentclass[12pt,reqno]{amsart}

\usepackage{latexsym}
\usepackage{amssymb}

\renewcommand{\Re}{{\operatorname{Re}\,}}
\renewcommand{\Im}{{\operatorname{Im}\,}}

\newcommand{\s}{{\mathbf s}}

\renewcommand{\epsilon}{\varepsilon}
\newcommand{\dist}{{\operatorname{dist}}}
\newcommand{\var}{{\operatorname{Var}}}
\newcommand{\Var}{{\bf Var}}

\newcommand{\sm}{\smallsetminus}
\newcommand{\szego}{Szeg\H{o} }

\newcommand{\inv}{^{-1}}
\newcommand{\kahler}{K\"ahler }
\newcommand{\sqrtn}{\sqrt{N}}
\newcommand{\wt}{\widetilde}
\newcommand{\wh}{\widehat}
\newcommand{\PP}{{\mathbb P}}
\newcommand{\N}{{\mathbb N}}
\newcommand{\R}{{\mathbb R}}
\newcommand{\C}{{\mathbb C}}

\newcommand{\Z}{{\mathbb Z}}

\newcommand{\CP}{\C\PP}
\renewcommand{\d}{\partial}
\newcommand{\dbar}{\bar\partial}
\newcommand{\ddbar}{\partial\dbar}

\newcommand{\K}{{\mathbf K}}
\newcommand{\E}{{\mathbf E}}
\newcommand{\Heis}{{\mathbf H}}

\newcommand{\half}{{\textstyle \frac 12}}
\newcommand{\vol}{{\operatorname{Vol}}}

\newcommand{\codim}{{\operatorname{codim\,}}}
\newcommand{\SU}{{\operatorname{SU}}}
\newcommand{\FS}{{{\operatorname{FS}}}}
\newcommand{\supp}{{\operatorname{Supp\,}}}
\renewcommand{\phi}{\varphi}
\newcommand{\eqd}{\buildrel {\operatorname{def}}\over =}

\newcommand{\ccal}{\mathcal{C}}
\newcommand{\dcal}{\mathcal{D}}
\newcommand{\ecal}{\mathcal{E}}

\newcommand{\hcal}{\mathcal{H}}

\newcommand{\lcal}{\mathcal{L}}

\newcommand{\ncal}{\mathcal{N}}
\newcommand{\ocal}{\mathcal{O}}

\newcommand{\scal}{\mathcal{S}}

\newcommand{\ycal}{\mathcal{Y}}

\newcommand{\al}{\alpha}
\newcommand{\be}{\beta}
\newcommand{\ga}{\gamma}

\newcommand{\La}{\Lambda}
\newcommand{\la}{\lambda}
\newcommand{\ep}{\varepsilon}
\newcommand{\de}{\delta}
\newcommand{\De}{\Delta}
\newcommand{\om}{\omega}
\newcommand{\Om}{\Omega}
\newcommand{\di}{\displaystyle}
\newtheorem{theo}{{\sc Theorem}}[section]

\newtheorem{cor}[theo]{{\sc Corollary}}

\newtheorem{lem}[theo]{{\sc Lemma}}
\newtheorem{prop}[theo]{{\sc Proposition}}

\newenvironment{rem}{\medskip\noindent{\it Remark:\/} }{\medskip}
\newtheorem{defin}[theo]{{\sc Definition}}

\title
{Number variance of random zeros on complex manifolds}

\author{Bernard Shiffman}
\author{Steve Zelditch}

\address{Department of Mathematics, Johns Hopkins University, Baltimore, MD
21218, USA} \email{shiffman@math.jhu.edu, zelditch@math.jhu.edu}

\thanks{Research of the first author partially supported by NSF grants
DMS-0100474 and DMS-0600982; research of the second author
partially supported by NSF grants  DMS-0302518 and DMS-0603850.}

\begin{document}

\begin{abstract}

We show that  the variance of the
number of simultaneous zeros of $m$ i.i.d.\ Gaussian random polynomials
of degree $N$ in an open set $U \subset \C^m$ with smooth boundary is asymptotic to
$N^{m-1/2} \,\nu_{mm}
\,\vol(\d U)$, where $\nu_{mm}$ is a universal constant depending
only on the dimension $m$. We also give formulas for the variance of the volume
of the set  of simultaneous
zeros in $U$ of $k < m$ random degree-$N$ polynomials on $\C^m$.  Our results hold
 more generally for the simultaneous  zeros of random
holomorphic sections of the $N$-th power of any positive line bundle
over any $m$-dimensional compact \kahler manifold.

\end{abstract}

\maketitle

\tableofcontents

\section{Introduction}

This article is concerned with the asymptotic statistics of the
number $\ncal_N^U(p_1^N,\dots, p_m^N)$  of  zeros in an open
set
$U \subset \C^m$ of a full system $\{p^N_j\}$ of $m$ Gaussian
random polynomials (or more generally,  of sections of a
holomorphic line bundle over a \kahler manifold $M_m$) as the
degree $N \to \infty$. In earlier work \cite{SZ}, we proved that
the zeros become uniformly distributed in $U$ with respect to the
natural volume form. The main result of this article (Theorem
\ref{number}) gives an asymptotic formula for the
 variance of
$\ncal_N^U(p_1^N,\dots,p_m^N)$ for open sets with piecewise smooth
boundary.  We also give analogous  results for the volume of the
simultaneous zero set of $k<m$ polynomials or sections (Theorem
\ref{volume}). Our results show that the zeros of a random system
are close to the expected distribution, i.e.  number statistics
are `self-averaging',  and moreover the degree of self-averaging
increases with the dimension.

To introduce our results, let us start with the case of
polynomials in $m$ variables.   By homogenizing, we may identify
the space of polynomials of degree $N$  in $m$ variables with the
space $H^0(\CP^m, \ocal(N))$ of holomorphic sections of the $N$-th
power of the hyperplane section bundle over $\CP^m$. This space
carries a natural $\SU(m+1)$-invariant inner product and
associated Gaussian measure $\gamma_N$. To each $m$-tuple of
degree $N$ polynomials $(p_1^N,\dots,p_m^N)$, we associate its
zero set $\{ p_1^N(z)=\cdots=p_m^N(z)=0\}$, which is almost always
discrete, and thus obtain a random point process on $\CP^m$. We
denote by $$Z_{p_1^N,\dots,p_m^N}:= \{z\in M: p_1^N(z)= \cdots
=p_m^N(z)=0\}\;$$ the set of zeros and by \begin{equation}
\label{ZCURRENT} \langle \left[Z_{p_1^N,\dots,p_m^N}\right], \psi
\rangle: = \sum_{z \in Z_{p_1^N,\dots,p_m^N}} \psi(z), \;\;\; \psi
\in C(M)
\end{equation} the sum of point masses at the joint zeros. It
easily follows from  the $\SU(m+1)$-invariance of $\ga_N$ that the
expected value  of this measure  is a multiple of the
Fubini-Study volume form, i.e.
\begin{equation}\label{MEASMEAN}\E \left[ Z_{p_1^N,\dots,p_m^N}\right] = N^m  \left(\frac
1\pi\,\om_\FS\right)^m\;,\end{equation}  where $\om_\FS= \frac i2
\ddbar \log |z|^2$ is the Fubini-Study metric on $\CP^m$. Here,
$\E(X)$ denotes the expected value of a random variable $X$.

 Given a measurable set $U$, the
random variable counting the number of zeros of the polynomial
system in $U$ is defined by
\begin{equation*} \ncal_N^U(p_1^N,\dots,p_m^N) := \#\{z\in U:
p_1^N(z)=\cdots=p_m^N(z)=0\}.
\end{equation*}
 Clearly, $\ncal^U_N$  is discontinuous along the
set of polynomials having a zero on the boundary $\partial U$.
Integrating (\ref{MEASMEAN}) over $U$ gives
\begin{equation}\label{MEAN}\E (\ncal^U_N) = N^m\int_U \left(\frac
1\pi\,\om_\FS\right)^m\;.\end{equation}

Formula \eqref{MEAN} has a counterpart for   Gaussian random
holomorphic sections of powers of any Hermitian holomorphic line
bundle $(L,h)$ with positive curvature $\Theta_h$ over any
$m$-dimensional compact \kahler manifold $M$ with \kahler form
$\om= \frac i2 \Theta_h$.  The Hermitian inner product
\begin{equation}\label{inner}\langle
\sigma_1, \bar \sigma_2 \rangle = \int_M h^N(\sigma_1, \bar
\sigma_2)\,\frac 1{m!}\om^m\;,\qquad \sigma_1,  \sigma_2 \in
H^0(M,L^N)\;,\end{equation} induces the complex Gaussian
probability measure
\begin{equation}\label{gaussian}d\gamma_N(s^N)=\frac{1}{\pi^m}e^
{-|c|^2}dc\,,\qquad s^N=\sum_{j=1}^{d_N}c_jS^N_j\,,\end{equation}
on the space $ H^0(M, L^N)$  of holomorphic sections of $L^N$,
where $\{S_1^N,\dots,S_{d_N}^N\}$ is an orthonormal basis for
$H^0(M,L^N)$, and $dc$ denotes $2d_N$-dimensional Lebesgue
measure.   The Gaussian measure $\ga_N$ given by
\eqref{inner}--\eqref{gaussian} is called the {\it Hermitian
Gaussian measure induced by $h$\/}.  The Gaussian ensembles
$(H^0(M,L^N), \ga_N)$ were used in \cite{SZ,SZ2,BSZ1,BSZ2}; for
the case of polynomials in one variable, they are equivalent to
the $\SU(2)$ ensembles studied in \cite{BBL, Han, NV} and
elsewhere.  In \cite{SZ}, we showed that  the expected value of
the corresponding random variable $\ncal^U_N$ (where $U\subset M$)
on the ensemble $(H^0(M, L^N)^m,\ga_N^m)$  has the asymptotics
\begin{equation}\label{MEANa}\frac 1
{N^m}\,\E (\ncal^U_N) = \int_U \left(\frac
i{2\pi}\,\Theta_h\right)^m +\ O\left(\frac
1N\right)\;.\end{equation} Thus, zeros of Gaussian random systems
of sections of $L^N$ become uniformly distributed with respect to
the curvature volume form $\frac 1{m!}\om^m$, as  $N\to\infty$.

Our main result  gives an asymptotic formula for the {\it number
variance}
$$\var(\ncal^U_N): = \E \big( \ncal^U_N - \E(\ncal^U_N) \big)^2$$ in this general
setting:

\begin{theo}\label{number}  Let $(L,h)$ be a
positive Hermitian holomorphic line bundle over a compact $m$-dimensional \kahler
manifold $M$. We give $H^0(M,L^N)$ the Hermitian Gaussian measure
induced by $h$ and the \kahler form $\om= \frac i2 \Theta_h$. Let
$U$ be a domain in $M$ with piecewise $\ccal^2$ boundary and no
cusps. Then for $m$ independent random sections $s_j^N\in
H^0(M,L^N)$, $1\le j\le m$, the variance of the random variable
$$\ncal^U_N(s_1^N,\dots,s_m^N):=\#\{z\in U:s_1^N(z)= \cdots
=s_m^N(z)=0\}$$ is given by
$$\var\big(\ncal^U_N\big) =
N^{m-1/2} \,\left[\nu_{mm} \,\vol_{2m-1}(\d U)
 +O(N^{-\frac 12 +\ep})\right]\;,$$ where $\nu_{mm}$ is a universal positive
constant. In particular,   $\nu_{11} =
\frac{\zeta(3/2)}{8\pi^{3/2}}$.
\end{theo}

Here, we say that $U$ has piecewise $\ccal^k$ boundary without
cusps if for each boundary point $z_0\in\d U$, there exists a (not
necessarily convex) closed polyhedral cone $K\subset\R^{2m}$ and a
$\ccal^k$ diffeomorphism $\rho:V\to \rho(V)\subset\R^{2m}$, where
$V$ is a neighborhood of $z_0$, such that $\rho (V\cap \overline
U) = \rho(V)\cap K$. By $O(N^{-\frac 12 +\ep})$, we mean a term
whose magnitude is less than $C_p N^p$ for some constant
$C_p\in\R^+$ (depending on $M,L,h,U$  as well as $p$), for all
$p>-\frac 12$.

As mentioned above, a special case of Theorem \ref{volume} gives
statistics for the number of zeros of systems of polynomials of
degree $N$.  Identifying polynomials on $\C^m$ of degree $N$ with
$H^0(\CP^m,\ocal(N))$ endowed with the Fubini-Study metric,  we
obtain an orthonormal basis of monomials (see \cite{BSZ2,SZ}):
$$\textstyle \left\{{N \choose
J}^{1/2}  z_1^{j_1}\cdots z_m^{j_m}\right\}_{|J|\le N}  \quad\big(
J=(j_1,\dots,j_m), \ |J|=j_1+\cdots+j_m,\ {N \choose J}= \frac
{N!}{(N-|J|)!j_1!\cdots j_m!}\,\big)$$

The polynomial case of Theorem \ref{number} then takes the form:
\begin{cor} \label{poly} Consider the Gaussian random polynomials
$$p^N_l(z_1,\dots,z_m) = \sum_{\{J\in\N^m:|J|\le N\}} c_{J,l}\, \textstyle  {N
\choose J}^{1/2} z_1^{j_1}\cdots z_m^{j_m}\qquad (1\le l\le
m)\;,$$ where the $c_{J,l}$ are independent complex Gaussian
random variables with mean 0 and variance 1. Let $U$ be a domain
in $\C^m$ with piecewise $\ccal^2$ boundary and no cusps. Then the
variance of the number $\ncal^U_N$ of zeros in $U$ of the degree-$N$
polynomial system $(p^N_1,\dots,p^N_m)$ is given by
$$\var\big(\ncal^U_N\big) =
N^{m-1/2} \,\left[\nu_{mm} \,\vol^{\CP^m}_{2m-1}(\d U)
 +O(N^{-\frac 12 +\ep})\right],$$
where $\vol^{\CP^m}_{2m-1}$ denotes the hypersurface
volume with respect to the Fubini-Study metric on $\CP^m$. \end{cor}

The variance $\var(\ncal^U_N)$ measures  the fluctuations in the
number of zeros in $U$ of random systems of polynomials or
sections. Theorem \ref{number} implies that the number of zeros in
$U$ is `self-averaging' in the sense that its fluctuations  are of
smaller order than its typical values. Recalling  \eqref{MEANa},
we have:

\begin{cor} \label{selfave}  Under the hypotheses of Theorem \ref{number} or
Corollary \ref{poly},
$$\frac{\left[\var\big(\ncal^U_N\big)\right]^{1/2}}{\E\big(\ncal^U_N\big)}
\sim  N^{-\frac m2 -\frac 14} \to 0 \qquad \mbox{as }\ N\to
\infty\;.$$
\end{cor}
By Corollary \ref{selfave} and the Borel-Cantelli Lemma applied to
the sets $\{ N^{-2m}(\ncal^U_N-\E\ncal^U_N)^2>\ep\}$, the zeros of
a {\it random sequence} $\{(s_1^N,\dots,s_m^N):N=1,2,3,\dots\}$ of
systems almost surely become uniformly distributed:
$$\frac 1{N^m}\;\ncal^U_N(s_1^N,\dots,s_m^N) \to m!\,\vol (U)
\qquad a.s.$$

Our proof of Theorem \ref{number} also yields asymptotic formulas
for the variance of volumes of simultaneous zero sets
$Z_{s_1^N,\dots,s_k^N}\cap U$, where
$$Z_{s_1^N,\dots,s_k^N}:= \{z\in M: s_1^N(z)= \cdots =s_k^N(z)=0\}\;.$$
of $k< m$ sections $s_1^N,\dots,s_k^N$. For $N$ sufficiently large
(so that $L^N$ is base point free), $Z_{s_1^N,\dots,s_k^N}$ is
almost always a complex codimension $k$ submanifold of $M$ (by
Bertini's Theorem) and
\begin{equation}\label{vol}\vol_{2m-2k}(Z_{s_1^N,\dots,s_k^N}\cap U)=
\int_{Z_{s_1^N,\dots,s_k^N}\cap U} \frac
1{(m-k)!}\,\om^{m-k}\;.\end{equation} We have the following
asymptotic formula for the variance of the  volume in \eqref{vol}:

\begin{theo}\label{volume}  Let $1\le k\le m$.  With the same notation
and hypotheses as in Theorem \ref{number}, for $k$ independent
random sections $s_j^N\in H^0(M,L^N)$, $1\le j\le k$, we have
$$\var\big(\vol_{2m-2k}[Z_{s_1^N,\dots,s_k^N}\cap U]\big) =
N^{2k-m-1/2} \,\left[\nu_{mk} \,\vol_{2m-1}(\d U)
 +O(N^{-\frac 12 +\ep})\right]\;,$$ where $\nu_{mk}$ is a universal positive  constant;  
in particular, 
$$\nu_{m1} =
\frac{\pi^{m-5/2}}{8}\;\zeta(m+\textstyle \half)\;.$$
\end{theo}

Theorem \ref{number} is the case $k=m$ of Theorem \ref{volume},
where the 0-dimensional volume is the counting measure:
$$\vol_0(Z_{s_1^N,\dots,s_m^N}\cap U) = \ncal^U_N(s_1^N,\dots,s_m^N)\;.$$  
  In fact,  Theorem
\ref{volume} is our means of proving Theorem \ref{number}, as our analysis
makes use of induction on the codimension.

The universal constants $\nu_{mk}$ in Theorem \ref{volume} must be
nonnegative, since variances (of nonconstant random variables) are
positive. In \S\ref{positivity}, we give explicit formulas for the $\nu_{mk}$
and show that they are strictly positive for all
$m$ and $k$.  

In the remainder of the introduction, we discuss  related  results on
 the variance problems studied in this article and indicate some
 key ideas in the proofs. In particular, we indicate which aspects
 of the results and methods are essentially complex analytical and
 which aspects are mainly probabilistic.

 The first results on number variance in domains
 appear to be due to   Forrester and Honner \cite{FH} for
 certain one-dimensional Gaussian ensembles of random complex polynomials. 
They gave an intuitive derivation of
the leading term of the asymptotic
 formula in dimension one,  which is proved in Theorem \ref{number}.
 Peres-Virag have precise results on numbers of zeros of Gaussian
 random analytic functions on the unit
 disc for a certain ensemble with a determinantal zero point
 process \cite{PV}.  To our knowledge, there are no prior results
 on number or volume variance in higher dimensions.

Variance asymptotics have also been studied for smooth analogues
of numbers statistics, namely for the random variables
(\ref{ZCURRENT}) with $\psi \in C^{\infty}(M)$. A rather simple
and non-sharp estimate on the variance for the smooth linear
statistics was given in our article \cite{SZ} to show that the
codimension-one  zeros of a random sequence $\{\s^N\}$ almost
 surely become uniformly distributed.
Smooth linear statistics were then studied in depth for certain
model one-dimensional Gaussian analytic functions by
Sodin-Tsirelson \cite{ST}) as a key  ingredient in their proof of
asymptotic normality for linear statistics. For their model
ensembles, they gave a sharp estimate for the variance of
$(Z_{s^N},\phi)$ and determined the leading term. (The constant
$\frac{\zeta(3)}{16\pi}$ was given for model ensembles in a
private communication from M. Sodin.)

We now  discuss some key ideas in the proofs,  and also their
relation to Sodin-Tsirelson \cite{ST} and to our prior work
\cite{BSZ1, SZa}. Apart from a rather routine computation (Lemma
\ref{varintb}) for the expected value $\E\big(\log |Y_1|\, \log
|Y_2|\big)$ where $Y_j$ are complex normal random variables, the
principal ingredients in our work are purely complex analytical.
The key ones are:

\begin{itemize}

\item A bipotential formula  for the variance current of one
section  (Theorem \ref{BIPOT}) or several sections (Theorem
\ref{variant}), and the closely related formulas for the pair
correlation current $\K_{21}^N$ (cf.\ (\ref{pcc})) in the
codimension one case (Proposition \ref{KN21});

\item Analysis of the singularities along the diagonal of
the pair correlation and variance currents, particularly in the point case
(maximal codimension)  where these currents contain a delta-function
along the diagonal in $M
\times M$ (Theorem \ref{variant1}). This analysis is necessary to verify the
formula
$\K_{2k}^N =
\left[\K_{21}^N\right]^{\wedge k}$   (see (\ref{powers})), where $\K_{2k}^N$ is the
pair correlation current   for the simultaneous zeros of $k$ random sections or
polynomials of degree $N$, and to define the product
$\left[\K_{21}^N\right]^{\wedge k}$.

\item Application of the rapid (in fact, exponentially fast)
off-diagonal asymptotics of the \szego kernel of \cite{SZ2} as the
degree $N \to \infty$ to obtain asymptotics of the variance
current and number variance.

\end{itemize}

Let us discuss these items in more detail. The   `bipotential' for
the pair correlation current was introduced in \cite{BSZ1}; in the
notation used here, the bipotential is a function $ Q_N(z, w)$
such that
\begin{equation}\label{biST}\Delta_z \Delta_w  Q_N (z,w) =
K_{21}^N(z,w)\,, \end{equation} where $K_{21}^N$ is the `pair
correlation function' (see \eqref{pcc}--\eqref{pcf}) for the zeros
of degree $N$ polynomials on $\CP^m$ or of holomorphic sections of
$L^N\to M$ . In this article, we show that $Q_N(z,w)$ is actually
a {\it pluri-bipotential\/} for the {\it variance current\/} $\Var
(Z_{s^N})$ of
$Z_{s^N}$ (see Theorem \ref{BIPOT}), so we can use $Q_N(z,w)$ to give an
explicit formula (Theorems
\ref{varint2}) for the variance of the zeros of $k\le m$
independent sections of
$H^0(M,L^N)$, for any codimension $k$.

The existence of the bipotential makes  essential use of the
complex analyticity of the polynomials and sections, in particular
the Poincar\'e-Lelong formula.  It is not clear if there exists a
useful generalization of (\ref{biST}) to the non-holomorphic
setting. Moreover, the analysis of the singularities of the
variance current requires a detailed study of intersections of
complex analytic  varieties and currents of fixed degree. Part of
the length of this paper is due to the lack of a prior  reference
for the relevant facts on smoothing and intersection of currents.
We hope that the analysis in \S \ref{s-bi} will be useful in other
applications.

The bipotential for the variance was also used implicitly by
Sodin-Tsirelson \cite{ST}, where it is defined as a power series
in the \szego kernel for $\ocal(N)\to\CP^1$. They used it to
obtain the first sharp formula for the variance of certain model
one-dimensional  random analytic functions, and further used it to
prove asymptotic normality of smooth linear statistics.

 From the bipotential  formulas, and the analysis of the singularities
 of the variance current,   the asymptotics of the variance
 are reduced to the
off-diagonal asymptotics for the \szego kernel $\Pi_N(z,w)$  (two
point function). Our main results are proved by applying the
off-diagonal asymptotics of $\Pi_N(z,w)$  in \cite{SZ2} to obtain
asymptotics of $Q_N(z,w)$ and then of the variance. One consequence of these
asymptotics, which is of independent interest, is that the normalized pair
correlation function $\wt K_{2m}^N$ rapidly approaches 1 off the diagonal as
$N$ increases:
\begin{equation}\label{cordecay0}\wt K_{2m}^N(z,w)  = 1 + O\left({N^{-\infty}}
\right)\;,\quad
\mbox{for }\ \dist(z,w)\ge N^{-1/2+\ep}\;.\end{equation}
(See Corollary \ref{cordecay} for a precise statement and the definition of $\wt
K_{2m}^N$.) These
asymptotics make essential use of the holomorphic setting, and
moreover of the positive curvature of the line bundles involved.
The two-point function $\Pi_N(z,w)$ and consequently the
correlation functions may behave quite differently in
non-holomorphic cases (the  two point kernel can decay at   only a
power law rate in real cases) or even for random
polynomials of degree $N$ on domains in $\C^m$ with the flat metric, i.e. with an inner
product independent of $N$. Subsequent to \cite{SZ2}, sharper
off-diagonal estimates for $\Pi_N(z,w)$ with exponentially small
remainder estimates away from the diagonal (i.e., when dist$(z,w)
\gg \frac 1 \sqrtn$) were also given in \cite{DLM, MaMa}, and they
would improve the decay of correlations in \eqref{cordecay0};
we state the result as above, since the estimates of \cite{SZ2}
already suffice for our applications.

The paper is organized as follows: In \S \ref{background}, we
review the formulas for the expected zero currents and describe
the asymptotics of the \szego kernel for powers of a line bundle.
In \S \ref{s-bi}, we define the variance current (in codimension
one) and introduce the pluri-bipotential $Q_N(z,w)$ for the
variance current and study its off-diagonal asymptotics. Next, we
provide our explicit formula  for the variance (Theorem \ref{varint2}). In
\S \ref{s-number}, this formula and the asymptotics of the
pluri-bipotential are applied to prove  Theorems
\ref{number} and \ref{volume} on number and volume variance.
Finally, in the Appendix (\S \ref{s-proof}), we review and to some
degree sharpen the derivation of the off-diagonal asymptotics in
\S \ref{off}.

We end the introduction with a word on the relation of this article to its
predecessors posted on arXiv.org.  The first predecessor of this
article is our preprint \cite{SZa}, in which we proved the
codimension $k=1$ case of Theorem \ref{volume}. This prior article
did not contain results on the point case in higher dimensions
since, as we wrote there, ``new technical ideas seem to be
necessary to obtain limit formula for the intersections of the
random zero currents $Z_{s_j}$."  The present article furnishes the
necessary new methods (cf.\ \S
\ref{s-bi}).  The preprint
\cite {SZa} also
extended the Sodin-Tsirelson asymptotic normality result
for smooth statistics \cite{ST}  to general one-dimensional ensembles and to
codimension one zero sets in higher dimensions.   It remains an
interesting open problem to generalize asymptotic normality to the point
case in higher dimension.  The original arXiv.org posting of the present
article (arxiv.org/abs/math/0608743v1) also contains results on smooth linear
statistics  and on random entire functions on
$\C^m$ and on certain noncompact complete \kahler
manifolds.   The results from these prior preprints on asymptotic
normality in codimension one and on smooth linear statistics are given in
our article \cite{SZ4}; the results on the noncompact case will be presented
elsewhere. 

\medskip\noindent{\it Acknowledgment:\/}  We thank M. Sodin for
discussions of his work with B. Tsirelson.   We also thank John Baber for
computing $\nu_{m2}$ using formulas from the
previous version of this article (arxiv.org/abs/math/0608743v2), where
positivity of the $\nu_{mk}$ was not established for all codimensions.  Baber's
computation suggested that the formulas for
$\nu_{mk}$ could be simplified, and we were then able to revise our
computations to show positivity and to obtain explicit formulas for
the constants in all codimensions.

 \section{Background}\label{background}

In this section we review the basic facts about the distribution
of zeros and the asymptotic properties of \szego kernels.

\subsection{Expected distribution of zeros}\label{EDZ}

In this section, we review the formula for the expected simultaneous zero
current of $k\le m$ independent Gaussian random  sections of the tensor
powers $L^N=L^{\otimes N}$ of a holomorphic line bundle
$L$ over an $m$-dimensional complex manifold $M$ (Corollary
\ref{indep}).  In order to
give a simple proof of our formula by induction on the codimension
$k$ of the zero set, we state our result in a more general form
(Proposition \ref{indep0}).  In the point case $k=m$, this formula was given by
Edelman-Kostlan \cite[Theorem~8.1]{EK} (for the essentially equivalent case
of trivial line bundles) using integral-geometric methods.

Throughout this paper, we let  $(L,h)$ be a Hermitian holomorphic
line bundle over a compact complex manifold $M$. We let $\scal$ be
a  subspace of $H^0(M,L)$, endowed with an (arbitrary)
Hermitian inner product. The inner product induces the complex Gaussian
probability measure
\begin{equation}\label{gaussian0}d\gamma(s)=\frac{1}{\pi^m}e^
{-|c|^2}dc\,,\qquad s=\sum_{j=1}^{n}c_jS_j\,,\end{equation} on
$\scal$, where $\{S_j\}$ is an orthonormal basis for $\scal$ and
$dc$ is $2n$-dimensional Lebesgue measure. This Gaussian is
characterized by the property that the $2n$ real variables $\Re
c_j, \Im c_j$ ($j=1,\dots,n$) are independent Gaussian random
variables with mean 0 and variance $\half$; i.e.,
$$\E c_j = 0,\quad \E c_j c_k = 0,\quad  \E c_j \bar c_k =
\de_{jk}\,.$$

 We let
\begin{equation}\label{sdef0}\Pi_\scal(z,z) =\E_\ga\left(
\|s(z)\|_h^2\right)=
\sum_{j=1}^n
\|S_j(z)\|_h^2\;,\qquad z\in M\;,\end{equation} denote the `\szego kernel'
for $\scal$ on the diagonal.  We now consider a local holomorphic frame
$e_L$ over a trivializing chart
$U$, and we write $S_j = f_j e_L$ over $U$. Any section
$s\in\scal$ may then be written as
\begin{equation}\label{frame}s = \langle c, F \rangle e_L\;, \quad
\mbox{where\ \ \ } F=(f_1,\dots,f_k)\;,\quad\langle c,F \rangle = \sum_{j =
1}^n c_j f_j\;.\end{equation} If  $s = f e_L$,  its Hermitian norm is given
by
$\|s(z)\|_h = a(z)^{-\half}|f(z)|$ where \begin{equation}
\label{a} a(z) = \|e_L(z)\|_h^{-2}\;. \end{equation} Recall that
the curvature form of $(L,h)$ is given locally by
$$\Theta_h= \ddbar \log a\;,$$ and the
{\it Chern form\/} $c_1(L,h)$ is given by
\begin{equation}\label{chern}c_1(L,h)=\frac{\sqrt{-1}}{2 \pi}
\Theta_h=\frac{\sqrt{-1}}{2 \pi}\d\dbar\log a\;.\end{equation}

Using  standard notation, we let $\ecal^{p,q}(X)$ and $\dcal^{p,q} (X)$ denote the
spaces of
${\ccal}^\infty$
$(p,q)$-forms and  compactly supported ${\ccal}^\infty$ $(p,q)$-forms, respectively,
on a complex manifold $X$, and we let
$\dcal'{}^{p,q}(X) =
\dcal^{m-p,m-q}(X)'$ denote the space of $(p,q)$-currents on $X$; $(T,\phi) =
T(\phi)$ denotes the pairing of
$T\in \dcal'{}^{p,q}(X)$ and $\phi \in
\dcal^{m-p,m-q}(X)$. If $Y$ is a complex submanifold of $X$ of codimension $p$, we
let $[Y]\in\dcal'^{p,p}(X)$ denote the current of integration over $Y$ given by
$([Y],\phi) = \int_Y\phi$.  (For a further description of currents on complex
manifolds, see
\cite[Ch.~3]{GH}.)

We now suppose that $\scal$ is {\it base point free\/}; i.e., the set $\{z\in
M: s(z)=0, \ \forall s\in\scal\}$ is empty.  By Bertini's theorem (or by an
application of Sard's theorem), for almost all $k$-tuples
$(s_1,\dots,s_k)\in\scal^k$, the simultaneous zero set
$\{z\in M:s_1(z)=\cdots=s_k(z)=0\}$ is a complex
submanifold of
$M$ of codimension
$k$, and we let $Z_{s_1,\dots,s_k}:= \big[\{s_1(z)=\cdots=s_k(z)=0\}\big]$
denote the current of integration over the zero
set:
\begin{equation*}(Z_{s_1,\dots,s_k},\phi)= \int_{\{ s_1(z)=\cdots =
s_k(z)=0\}} \phi(z), \qquad \phi\in
\dcal^{m-k,m-k}(M)\;. \end{equation*} The current $Z_{s_1,\dots,s_k}$ is
well-defined for almost all
$s_1,\dots,s_k$, and we regard $Z_{s_1,\dots,s_k}$ as a current-valued random
variable.   For the point case $k=m$,
$Z_{s_1,\dots,s_m}$ is a measure-valued random variable:
\begin{equation*}(Z_{s_1,\dots,s_m},\phi)= \sum_{ s_1(z)=\cdots =
s_m(z)=0} \phi(z),\qquad \phi\in\dcal(M)\;.  \end{equation*}

The
current of integration $Z_s$ over the zeros of one section $s\in\scal$,
written locally as $s=fe_L$, is then given by the {\it Poincar\'e-Lelong
formula\/}  (see
\cite[p.~388]{GH}):
\begin{equation} Z_{s} =
\frac{\sqrt{-1}}{ \pi } \partial \bar{\partial}\log |f| = \frac{\sqrt{-1}}{
\pi } \partial
\bar{\partial}\log  \left\|s\right\|_{h} + c_1(L,h)\;,
\label{Zs}
\end{equation} where the second equality is a consequence of
\eqref{a}--\eqref{chern}.

We begin with a general form of the {\it probabilistic Poincar\'e-Lelong formula\/}
from
\cite{SZ} (see also
\cite{BSZ1,BSZ2}) for the expected  value of a random
zero divisor:

\begin{prop}\label{EZ}Let $(L,h)$ be a Hermitian line bundle on a
compact \kahler manifold $M$. Let $\scal$ be a base-point-free
subspace of $ H^0(M,L)$ endowed with a Hermitian inner
product and we let $\ga$ be the induced  Gaussian probability measure
on $\scal$. Then the expected zero current of a random section
$s\in\scal$ is given by
\begin{eqnarray*}\E_\ga(Z_s)  &=&\frac{\sqrt{-1}}{2\pi}
\partial
\bar{\partial} \log \Pi_{\scal}(z, z)+c_1(L,h)\;.\end{eqnarray*}
\end{prop}

The formula of the proposition was essentially given in
\cite[Prop.~3.1]{SZ}.  For completeness, we give a short proof: It suffices to
verify the identity over a trivializing neighborhood $U$.  As above, we let
$\{S_j\}$ be an orthonormal basis for $\scal$,  we  write
$S_j = f_j e_L$ (over $U$), and we let $F=(f_1,\dots,f_n)$.  As in
\cite{SZ},  we then write $F(z)= |F(z)| u(z)$ so that $|u|
\equiv 1$ and
\begin{equation}\label{2terms}\log  | \langle c, F \rangle| = \log |F| +
\log  | \langle c, u \rangle|\;.\end{equation}

Thus by (\ref{Zs}), we have
\begin{eqnarray*}\big(\E_\ga(Z_s),\phi\big)&=&\frac{\sqrt{-1}}{ \pi}
\int_\scal \left(\log  | \langle c,F \rangle| , \ddbar\phi\right)
d\gamma\\ &=&\frac{\sqrt{-1}}{ \pi}  \left(\log  |F| ,
\ddbar\phi\right)  +\frac{\sqrt{-1}}{ \pi}\int_\scal\left(
\log  | \langle c, u
\rangle|, \ddbar\phi\right) d\gamma \;,
\end{eqnarray*}
for all test forms $\phi\in\dcal^{m-1,m-1}(U)$.

A key point is that $\langle c,u(z)\rangle$ is a standard (mean 0, variance 1)
complex Gaussian random variable for all $z\in
U$ (since $u(z)$ is  unit vector) and hence $\E\big(\log  |
\langle c, u(z)
\rangle|\big)$ is a universal constant $C$ independent of $z$.

Thus
\begin{eqnarray*}\int_\scal\left(
\log  | \langle c, u
\rangle|, \ddbar\phi\right) d\gamma &=& \int_\scal d\gamma(c)
\int_M \log  | \langle c, u
\rangle|\, \ddbar\phi\\&=& \int_M\left[ \int_\scal
 \log  | \langle c, u
\rangle|d\gamma(c)\right] \ddbar\phi\\&=&C \int_M  \ddbar\phi\ =
\ 0\;.
\end{eqnarray*}
Fubini's Theorem can be applied above
since
$$\int_{M\times \scal}\left|\log  | \langle c, u
\rangle|\, \ddbar\phi\right|\;d\gamma_N(c)= \left(\int_\C
\big|\log  | \zeta|\big|\,\frac 1\pi\,e^{-|\zeta|^2} d\zeta\right)
\left(\int_M |\ddbar\phi|\right)< +\infty\;.$$

Therefore
\begin{equation}\label{Egamma}\E_\ga(Z_{s})=\frac{\sqrt{-1}}{2 \pi} \partial
\bar{\partial}\log  |F|^2  =\frac{\sqrt{-1}}{2 \pi} \partial
\bar{\partial}\left(\log \sum\|S_j\|_{h}^2+N\log
a\right)\;.\end{equation} The formula of the proposition then
follows from \eqref{sdef0}, \eqref{chern} and \eqref{Egamma}.
\qed

\begin{rem} Proposition \ref{EZ} also holds if $M$ is noncompact and without
the assumption   that $\scal$ is base point free or even finite
dimensional.
\end{rem}

Next, we give a general result on the expected value of the simultaneous
zero current of
$k$ independent random holomorphic sections:

\begin{prop} \label{indep0} Let $(L,h)\to M$, $\scal\subset H^0(M,L)$,
and $\ga$ be given as in Proposition \ref{EZ}, and let $1\le k\le m$.  Then
the expected value of the simultaneous zero current of  $k$ independent
random sections $s_1,\dots,s_k$ in $\scal$ is given by
$$\E_{\ga^k}\big(Z_{s_1,\dots,s_k}\big) = \left(\frac{\sqrt{-1}}{2\pi}
\partial
\bar{\partial} \log \Pi_{\scal}(z, z)+c_1(L,h)\right)^k\;.$$\end{prop}

The proposition is a formal consequence of Proposition \ref{EZ}
and the independence of the sections $s_j$, but needs a proof
since the wedge products of currents is not always defined. We
give here a simple induction proof without using the theory of
wedge products of singular currents.

\begin{proof} Let $\om$ be the  \kahler
form on $M$. We first note that
\begin{equation}
\label{ind1} \big(Z_{s_1,\dots,s_k},\om^{m-k}\big) = \int_M
c_1(L,h)^k \wedge \om^{m-k}\;,\end{equation} whenever the
$Z_{s_j}$ are smooth and intersect transversely.  The identity
\eqref{ind1} is a consequence of the fact that the current
$Z_{s_1,\dots,s_k}$ and the smooth form $c_1(L,h)^k$ are in the
same de Rham cohomology class.  Equation \eqref{ind1} can also be
verified by induction: the case $k=1$ follows immediately from the
Poincar\'e-Lelong formula \eqref{Zs} and the fact that $\om$ is
closed; assuming the result for $k-1$ sections on $Z_{s_1}$, we
have
$$\big(Z_{s_1}\cap Z_{s_2,\dots,s_k},\om^{m-k}\big) = \int_{Z_{s_1}}
c_1(L,h)^{k-1} \wedge \om^{m-k} = \int_Mc_1(L,h)\wedge
c_1(L,h)^{k-1} \wedge \om^{m-k} \,,$$ which gives \eqref{ind1} for
$k$ sections.

Now let $\phi\in\dcal^{m-k,m-k}(M)$ be a test form. By
\eqref{ind1} and the formula for the volume of complex
submanifolds \eqref{vol}, we then have
\begin{eqnarray*} \left|\big(Z_{s_1,\dots,s_k},\phi\big)\right|
&=& \left|\int_{Z_{s_1,\dots,s_k}}\phi\right|\ \le\
\sup\|\phi\|\,\vol (Z_{s_1,\dots,s_k})\\ &=&
\frac{\sup\|\phi\|}{(m-k)!}\left(Z_{s_1,\dots,s_k},\om^{m-k}\right)\
=\ \frac{\sup\|\phi\|}{(m-k)!} \int_M c_1(L,h)^k \wedge \om^{m-k}
,\end{eqnarray*} for almost all $s_1,\dots,s_k$. Thus the random
variable $\big(Z_{s_1,\dots,s_k},\phi\big)$ is $L^\infty$, so its
expected value  is well defined.

We must show that \begin{equation}\label{ind3}
\E_{\ga^k}\big(Z_{s_1,\dots,s_k},\phi\big) = \int_M \al^k \wedge
\phi\;,\end{equation} where $$\al= \frac{\sqrt{-1}}{2\pi}
\partial
\bar{\partial} \log \Pi_{\scal}(z, z)+c_1(L,h)\;.$$
 We verify \eqref{ind3} by induction
on $k$: For $k=1$,  Proposition \ref{EZ}   yields  \eqref{ind3}. Let
$k\ge 2$ and suppose that \eqref{ind3} has been verified for $k-1$
sections. Choose $s_1\in\scal$ such that $Z_{s_1}$ is a
submanifold, and let $M'=Z_{s_1},\ s_j'=s_j|_{M'},\
\scal'=\scal|_{M'}$. We  give $\scal'$ the Gaussian measure
$\ga':=\rho_*\ga$, where $\rho:\scal\to\scal'$ is the restriction
map, and we note that
$$\Pi_{\scal'}(z',z') = \E_{\ga'}\left(
\|s'(z)\|_h^2\right)= \E_\ga\left( \|s(z)\|_h^2\right)=
\Pi_{\scal}(z',z')\;, \quad\mbox{for }\ z'\in M'\;.$$ By the
inductive assumption applied to $M',\scal'$, and noting that
$Z_{s_1,\dots,s_k}=Z_{s_2',\dots,s_{k}'}$, we have
\begin{equation}\label{ind2}\int_{\scal^{k-1}} \big(Z_{s_1,\dots,s_k},\phi\big)\,
d\ga(s_2)\cdots d\ga(s_{k})
=\E_{\ga'^{k-1}}\big(Z_{{s'_2},\dots,s'_{k}},\phi\big) =
\int_{Z_{s_1}} \al^{k-1}\wedge\phi\;.\end{equation} We average
\eqref{ind2} over $s_1$ and apply Proposition \ref{EZ} to conclude
that
$$\int_{\scal^{k}}
\big(Z_{s_1,\dots,s_k},\phi\big)\, d\ga(s_1)\cdots d\ga(s_{k}) =
\int_\scal \big(Z_{s_1},\al^{k-1}\wedge\phi\big)\, d\ga(s_1) =
\big(\al\,,\al^{k-1}\wedge\phi\big)= \int_M \al^k \wedge \phi,$$
which gives \eqref{ind3}.
\end{proof}

\subsubsection{Powers of a positive line bundle}\label{s-powers}  We now
specialize Proposition \ref{indep0} to our case of interest.  We suppose
that the line bundle $(L,h)$ is {\it positive\/}, i.e.\ the $(1,1)$-form
$c_1(L,h)$ is everywhere positive definite, and we give $M$ the \kahler
form $\om= \frac i2\Theta_h=\pi c_1(L,h)$. Recall that it is a consequence
of the Kodaira embedding theorem that for sufficiently large integers $N$,
the spaces of global sections of the tensor powers $L^N=L^{\otimes N}$ of
the line bundle are base point free.  (In fact, the global sections give
a projective embedding
\cite[\S 1.4]{GH}.  The Kodaira embedding theorem is also a consequence of
the Tian-Yau-Zelditch theorem
\cite{Cat,Ti,Z}; see
\eqref{TYZ}.)  We recall that the Hermitian metric $h$ on $L$ induces
Hermitian metrics
$h^N$ on
$L^N$, and we have $$c_1(L^N,h^N) = N\,c_1(L,h)= \frac N\pi \,\om\,.$$

We give $H^0(M,L^N)$ the Hermitian inner product induced by the metrics
$h,\om$, as defined by
\eqref{inner}; this inner product induces
the Hermitian Gaussian measure $\ga_N$  given by \eqref{gaussian}.
Considering the spaces $\scal_N=H^0(M,L^N)$, we have the \szego kernels (on
the diagonal)
\begin{equation}\label{sdef} \Pi_N(z,z):=
\Pi_{\scal_N}(z,z)=\sum_{j=1}^{d_N}
\|S^N_j(z)\|^2_{h^N} \;,\end{equation} where
$\{S_1^N,\dots,S_{d_N}^N\}$ is an orthonormal basis for
$H^0(M,L^N)$ with respect to the Hermitian inner product \eqref{inner}.
These
\szego kernels were analyzed in
\cite{Z,BSZ1,SZ, SZ2} by viewing them as orthogonal projectors on
$\lcal^2(X)$, where
$X\to M$ is the circle bundle of unit vectors of $L\inv$.  We give this
description of $\Pi_N$ in \S \ref{off} below.

Applying Proposition \ref{indep0} to the line bundles $L^N$ and the spaces
$H^0(M,L^N)$ of holomorphic sections, we obtain:

\begin{cor} \label{indep} Let $(L,h)\to (M,\om)$ be
as in Theorem
\ref{number}, and let $\ga_N$ be the Hermitian Gaussian measure on
$H^0(M,L^N)$. Then for $1\le k\le m$ and $N$ sufficiently large, we have
$$\E_{\ga_N^k}\big(Z_{s_1^N,\dots,s_k^N}\big) = \big(\E_{\ga_N}Z_{s^N}\big)^k =
 \left(   \frac i{2\pi}  \ddbar
\log \Pi_{N}(z,z) + \frac N\pi \omega\right)^k\,.$$\end{cor}

\subsection{Off-diagonal asymptotics for the \szego kernel}\label{off}

As in \cite{Z,SZ, BSZ1} and elsewhere,  we analyze  the \szego
kernel for $H^0(M,L^N)$ by lifting it to the circle bundle
$X{\buildrel {\pi}\over \to} M$  of unit vectors in the dual
bundle $L\inv\to M$ with respect to $h$. In the standard way
(loc.\ cit.), sections of $L^N$ lift to equivariant functions on
$X$. Then $s\in H^0(M,L^N)$ lifts to a $CR$ holomorphic functions
on $X$ satisfying  $\hat s(e^{i\theta}x)= e^{iN\theta}\hat s(x)$.
We denote the space of such functions by  $ \hcal^2_N(X)$. The
{\it \szego projector\/} is the orthogonal projector
$\Pi_N:\lcal^2(X)\to\hcal^2_N(X)$, which is given by the {\it
\szego kernel}
$$\Pi_N(x,y)=\sum_{j=1}^{d_N} \wh S^N_j(x)\overline{\wh S^N_j(y)}
\qquad (x,y\in X)\;.$$ (Here, the functions $\wh S^N_j$ are the
lifts to $\hcal^2_N(X)$ of the orthonormal sections $S_j^N$; they
provide an orthonormal basis for $\hcal^2_N(X)$.)

Further, the covariant derivative $\nabla s$  of a section $s$ lifts
to the horizontal derivative $\nabla_h \hat{s}$ of its
equivariant lift $\hat{s}$ to $X$; the horizontal derivative is of the
form
\begin{equation}\label{HORDER} \nabla_h \hat s
=\sum_{j=1}^m\left( \frac{\d \hat s}{\d z_j} -A_j\frac{\d\hat s}{\d
\theta}\right)dz_j.
\end{equation}
For further discussion and  details on lifting sections, we refer
to \cite{SZ}.

Our pluri-bipotential for the variance described in \S \ref{s-bi}
is based on the {\it normalized \szego kernels}
\begin{equation}\label{PN} P_N(z,w):=
\frac{|\Pi_N(z,w)|}{\Pi_N(z,z)^\frac 12 \Pi_N(w,w)^\frac
12}\;,\end{equation} where we write
$$|\Pi_N(z,w)|:=|\Pi_N(x,y)|\;,\quad z=\pi(x),\,w=\pi(y)\in M\;.$$
In particular, on the diagonal we have $\Pi_N(z,z)=\Pi_N(x,x)>0$.
Note that $\Pi_N(z,z)=\Pi_\scal(z,z)$ as defined in (\ref{sdef})
with $\scal=H^0(M,L^N)$.

In this section, we
use the off-diagonal asymptotics for $\Pi_N(x,y)$ from \cite{SZ2} to
provide the off-diagonal estimates for the  normalized \szego kernel
$P_N(z,w)$ that we need for our variance formulas.  Our estimates are of two
types: (1) `near-diagonal' asymptotics (Propositions
\ref{better}--\ref{best}) for
$P_N(z,w)$ where the distance
$\dist(z,w)$ between
$z$ and
$w$ satisfies an upper bound $ \dist(z,w)\le b\left(\frac
{\log N}{N}\right)^{1/2}$ ($b\in\R^+$); (2) `far-off-diagonal' asymptotics
(Proposition \ref{DPdecay}) where  $ \dist(z,w)\ge b\left(\frac
{\log N}{N}\right)^{1/2}$.

To describe the scaling asymptotics for the \szego kernel at a point $z_0
\in M$,  we choose a neighborhood
$U$ of $z_0$, a local normal coordinate chart
$\rho:U,z_0\to\C^m,0$ centered at $z_0$, and a {\it preferred\/}
local frame  at $z_0$, which we defined in \cite{SZ2} to be a
local frame $e_L$ such that
\begin{equation}\label{preferred}\|e_L(z)\|_h=1-\half \|\rho(z)\|^2
+ \cdots\;.\end{equation} For $u=(u_1,\dots,u_m)\in \rho(U),\
\theta\in (-\pi,\pi)$, we let
\begin{equation}\tilde\rho(u_1,\dots,u_m,\theta)=\frac
{e^{i\theta}} {|e^*_L(\rho\inv(u))|_h}
 e^*_L(\rho\inv(u))\in
X\,,\label{coordinates}\end{equation} so that
$(u_1,\dots,u_m,\theta)\in\C^m\times \R$ give local coordinates on
$X$. As in \cite{SZ2}, we write
$$\Pi_N^{z_0}(u,\theta; v,\phi)
=\Pi_N(\tilde\rho(u,\theta),\tilde\rho(v,\phi))\;.$$ Note that
$\Pi_N^{z_0}$ depends on the choice of coordinates and frame; we
shall assume that we are given normal coordinates and local frames
for each point $z_0\in M$ and that these normal coordinates and
local frames are smooth functions of $z_0$. The scaling
asymptotics of $\Pi_N^{z_0}(u,\theta; v,\phi)$ lead to the model
Heisenberg \szego kernel
\begin{equation}\label{heisen-N}\Pi^\Heis_N(z,\theta;w,\phi) =
e^{iN(\theta-\phi)}\sum_{k\in \N^m}S_k(z)\overline{S_k(w)} =
\frac{N^m}{\pi^m} e^{i  N(\theta - \phi)+N z \cdot \bar{w}- \frac
N2 (|z|^2+ |w|^2)} \end{equation} of level $N$ for the
Bargmann-Fock space of functions on $\C^m$ (see \cite{BSZ2}).

We shall apply the following (near and far) off-diagonal
asymptotics from \cite {SZ2}:

\begin{theo} \label{near-far} Let $(L,h)\to (M,\om)$ be as in Theorem
\ref{number}, and let $z_0\in M$.  Then using the above notation,
\begin{enumerate}
\item[i)] $\  N^{-m}\Pi_N^{z_0}(\frac{u}{\sqrtn},\frac{\theta}{N};
\frac{v}{\sqrtn},\frac \phi N)$
$$\begin{array}{l}
= \Pi^\Heis_1(u,\theta;v,\phi)\left[1+ \sum_{r = 1}^{k} N^{-r/2}
p_{r}(u,v) + N^{-(k +1)/2} R_{Nk}(u,v)\right]\;,\end{array}$$
where
 the $p_r$ are polynomials in $(u,v)$ of degree
$\le 5r$ (of the same parity as
$r$), and
$$|\nabla^jR_{Nk}(u,v)|\le C_{jk\ep b}N^{\ep}\quad \mbox{for }\
|u|+|v|<b\sqrt{\log N}\,,$$ for
$\ep,b\in\R^+$,  $j,k\ge 0$. Furthermore, the constant  $C_{jk\ep
b}$ can be chosen independently of $z_0$. \item[ii)] For
$b>\sqrt{j+2k+2m}\,$, $j,k\ge 0$, we have
$$ \left|\nabla^j_h
\Pi_N(z,w)\right|=O(N^{-k})\qquad \mbox{uniformly for }\ \dist(z,w)\ge
b\,\sqrt{\frac {\log N}{N}} \;.$$
\end{enumerate}\end{theo}

Here $\nabla^jR = \left\{\frac {\d^j R}{\d u^{K'} \d v^{K''}}:
|K'|+|K''|=j\right\}$, and $\nabla^j_h=(\nabla_h)^j$ denotes the
$j$-th iterated horizontal covariant derivative; see
(\ref{HORDER}). Theorem \ref{near-far} is equivalent to equations
(95)--(96) in \cite{SZ2}, where the result was shown to hold for
almost-complex symplectic manifolds. (The remainder in (i) was
given for $v=0$, but the proof holds without any change for $v\ne
0$.  Also the statement of the result was divided into the two
cases where the scaled distance is less or more, respectively,
than $N^{1/6}$ instead of $\sqrt{\log N}$ in the above
formulation, which is more useful for our purposes.)  A
description of the polynomials $p_r$ in part (i) is given in
\cite{SZ2}, but we only need the $k=0$ case in this paper. For the
benefit of the reader, we give a proof of Theorem \ref{near-far}
in \S\ref{s-proof}.

\begin{rem} The \szego kernel actually satisfies the sharper `Agmon decay
estimate' away from the diagonal:
\begin{equation}\label{agmon}\nabla^j\Pi_N(
z,\theta;w,\phi)= O\left(e^{-A_j\sqrtn\,\dist(z,w)} \right)\;,
\qquad j\ge 0\;.
\end{equation} In particular,
\begin{equation}\label{szegodecay} |\Pi_N( z,w)|=  O\left(e^{-A\sqrtn\,\dist(z,w)}
\right)\;.\end{equation} A short proof of \eqref{szegodecay} is
given in \cite[Th.~2.5]{Be}; similar estimates were established by
M.~Christ \cite{Ch}, H. Delin \cite{D}, and N. Lindholm\cite{Li}.
(See also \cite{DLM, MaMa} for off-diagonal exponential estimates
in a more general setting.)  We do not need Agmon estimates for
this paper; instead Theorem \ref{near-far} suffices.
\end{rem}

It follows from Theorem \ref{near-far}(i)) with $k=1$ that on the diagonal,
the \szego kernel is of the form
\begin{equation}\label{TYZ}\Pi_N(z,z)= \frac
1{\pi^m}N^m(1+O(N\inv))\,,\end{equation} which comprises the
leading terms of the Tian-Yau-Zelditch asymptotic expansion of the
\szego kernel \cite{Cat,Ti,Z}.  Applying \eqref{TYZ}, we obtain
the asymptotic formula from \cite{SZ} for the expected
simultaneous zero currents:

\begin{prop} \label{indepTYZ} {\rm \cite[Prop.~4.4]{SZ}} Let $(L,h)\to
(M,\om)$ be as in Theorem
\ref{number}, and let $1\le k\le m$.  Then for independent random sections
$s_1^N,\dots,s_k^N$ in $H^0(M,L^N)$, we have
$$\E\big(Z_{s_1^N,\dots,s_k^N}\big)= \frac {N^k}{\pi^k} \om^k
+O\left(N^{k-1}\right).$$\end{prop}

\begin{proof} By \eqref{TYZ},  $\ddbar\log \Pi_N(z,z) = O(N\inv)$.  The
asymptotics for $\E\big(Z_{s_1^N,\dots,s_k^N}\big)$ then follow from the
formula of Corollary
\ref{indep}.\end{proof}

We now state our far-off-diagonal decay estimate for $P_N(z,w)$, which
follows immediately from Theorem
\ref{near-far}(ii) and \eqref{TYZ}.

\begin{prop}\label{DPdecay} Let $(L,h)\to (M,\om)$ be as in Theorem
\ref{number}, and let $P_N(z,w)$ be the normalized \szego kernel for
$H^0(M,L^N)$ given by \eqref{PN}. For
$b>\sqrt{j+2k}$,
$j,k\ge 0$, we have
$$ \nabla^j
P_N(z,w)=O(N^{-k})\qquad \mbox{uniformly for }\ \dist(z,w)\ge
b\,\sqrt{\frac {\log N}{N}} \;.$$
\end{prop}
\medskip
The normalized \szego kernel $P_N$ also satisfies Gaussian
decay estimates valid very close to the diagonal. To give the estimate, we write by
abuse of notation,
\begin{equation}\label{aon}P_N(z_0+u,z_0+v):=
P_N(\rho\inv(u),\rho\inv(v))= \frac {|\Pi_N^{z_0}(u,0;
v,0)|}{\Pi_N^{z_0}(u,0; u,0)^{1/2}\Pi_N^{z_0}(v,0;
v,0)^{1/2}}\;.\end{equation} As an immediate consequence of Theorem
\ref{near-far}(i), we have:

\begin{prop} \label{better} Let $P_N(z,w)$ be as in Proposition
\ref{DPdecay}, and let $ z_0\in M$.  For $b,\ep>0,\  j\ge 0$, there is a constant
$C_j=  C_j({M,\ep,b})  $, independent of the point $z_0$, such that
\begin{eqnarray*}\textstyle  P_N\left(z_0+\frac u{\sqrtn},z_0+\frac v{\sqrtn}\right) &=&
e^{-\frac 12 |u-v|^2}[1 + R_N(u,v)]\\ &&\quad
|\nabla^jR_N(u,v)|\le C_j\,N^{-1/2+\ep}\quad \mbox{for }\
|u|+|v|<b\sqrt{\log N}\;.\end{eqnarray*}
\end{prop}

 As a corollary we have:

\begin{prop} \label{best}
The remainder $R_N$ in Proposition \ref{better} satisfies
$$ |R_N(u,v)|\le \frac {C_2}2\,|u-v|^2N^{-1/2+\ep}, \quad |\nabla R_N(u)| \le
C_2\,|u-v|\,N^{-1/2+\ep},
 \quad \mbox{for }\ |u|+|v|<b\sqrt{\log N}.$$\end{prop}

\begin{proof} Since $P_N\left(z_0+u,z_0+v\right) \le 1=
P_N\left(z_0+u,z_0+u\right)$, we conclude that $R_N(u,u)=0,$\break $
dR_N|_{(u,u)}=0$, and thus by  Proposition
\ref{better},
$$|\nabla R_N(u,v)| \le  \sup_{0\le t\le
1}|\nabla^2 R_N(u,(1-t)u+tv)|\,|u-v|\le C_2\,|u-v|\,N^{-1/2+\ep}\;.$$
 Similarly,
$$|R_N(u,v)| \le \half  \sup_{0\le t\le
1}|\nabla^2 R_N(u,(1-t)u+tv)|\,|u-v|^2\le \frac {C_2}2\,|u-v|^2\,N^{-1/2+\ep}\;.$$
\end{proof}

\section{A pluri-bipotential for the variance}\label{s-bi}

Our proof  of Theorems \ref{number} is based on a
pluri-bipotential given implicitly in \cite{SZ} for the variance
current for random zeros in codimension one. More generally, for
random codimension
$k$ zeros, we define the {\it variance current\/} of
$Z_{s_1^N,\dots,s_k^N}$ to be the current
\begin{equation}\label{vc}
\Var\big(Z_{s_1^N,\dots,s_k^N}\big): =  \E\big(Z_{s_1^N,\dots,s_k^N} \boxtimes
Z_{s_1^N,\dots,s_k^N}
\big) -
\E\big(Z_{s_1^N,\dots,s_k^N}\big)\boxtimes \E
\big(Z_{s_1^N,\dots,s_k^N}\big)\in \dcal'^{2k,2k}(M\times M).
\end{equation}
Here we  write
$$S\boxtimes T = \pi_1^*S \wedge \pi_2^*T \in \dcal'^{p+q}(M\times
M)\;, \qquad \mbox{for }\ S\in \dcal'^p(M),\ T\in \dcal'^q(M)\;,$$
where $\pi_1,\pi_2:M\times M\to M$ are the projections to the
first and second factors, respectively.  The variance for the
`smooth zero statistics' is given by:
\begin{equation}\label{var}\var\big(Z_{s_1^N,\dots,s_k^N},\phi\big)
=\left(\Var\big(Z_{s_1^N,\dots,s_k^N}\big),\;\phi\boxtimes
\phi\right)\;.\end{equation}  Conversely, \eqref{var} can be taken
as an equivalent definition of the variance current in terms of
$\var\big(Z_{s_1^N,\dots,s_k^N},\phi\big)$.

Theorem \ref{BIPOT} below gives a {\it pluri-bipotential\/} for
the variance current in codimension one, i.e.\  a function $Q_N\in
L^1(M\times M)$ such that
\begin{equation}\label{varcur}{\bf Var}\big(Z_{s^N}\big)=
(i\ddbar)_z\,(i\ddbar)_w \,Q_N(z,w) \;.\end{equation}

To describe our pluri-bipotential $Q_N(z,w)$, we define the function
\begin{equation}\label{Gtilde} \wt G(t):= -\frac 1{4\pi^2}
\int_0^{t^2} \frac{\log(1-s)}{s}\,ds\ =\ \frac 1{4\pi^2}
\sum_{n=1}^\infty\frac{t^{2n}}{n^2}\;,
\qquad 0\le t\le 1.\end{equation}   Alternatively,  
\begin{equation}\label{Gtilde1}\wt G(e^{-\la}) = -\frac
1{2\pi^2} \int_\la^\infty \log(1-e^{-2s})\,ds\;,\qquad
\la\ge 0\;.\end{equation} The function  $\wt G$ is a modification of
the function $G$ defined in \cite{BSZ1}; see \eqref{GGtilde}.

\begin{theo} \label{BIPOT} Let $(L,h)\to(M,\om)$ be as in Theorems
\ref{number}. Let $Q_N:M\times M\to [0,+\infty)$ be
the function given  by
\begin{equation}
\label{QN} Q_N(z,w)= \wt G(P_N(z,w)) = -\frac 1{4\pi^2}
\int_0^{P_N(z,w)^2} \frac{\log(1-s)}{s}\,ds\;, \end{equation}
where $P_N(z,w)$ is the normalized \szego kernel given by
\eqref{PN}. Then $${\bf Var}\big(Z_{s^N}\big)=
(i\ddbar)_z\,(i\ddbar)_w \,Q_N(z,w)\;.$$
\end{theo}
Theorem \ref{BIPOT} says that
\begin{equation} \label{varint1} \var( Z_{s^N},\phi) =
\big(-\d_z\dbar_z\d_w\dbar_w Q_N,\;\phi\boxtimes \phi\big)= \int_M
\int_M Q_N(z,w)\,(i\ddbar \phi(z))\, (i\ddbar \phi(w))
\;,\end{equation} for test forms
$\phi\in\dcal^{m-1,m-1}(M)$.
We note that $Q_N$ is $\ccal^\infty$ off the diagonal   for $N$
sufficiently large,  but is only
$\ccal^1$ and not $\ccal^2$ at all points on the diagonal in
$M\times M$, as the computations in \S \ref{Qasympt} show.
Additionally, its derivatives of order $\le 4$ are in
$L^{m-\ep}(M\times M)$ (see Lemma \ref{DISTM2}).

To begin the proof of the theorem, we write
\begin{equation}\label{Psi}\Psi_N=(S_1^N,\dots,S_{d_N}^N)\in
H^0(M,L^N)^{d_N}\;,\end{equation} where $\{S_j^N\}$ is an
orthonormal basis of $H^0(M,L^N)$. As in the proof of Proposition
\ref{EZ}, we write \begin{equation}\label{uN}\Psi_N(z) = |\Psi_N(z)|\,
u_N(z)\;,\end{equation} where
$|\Psi_N|:= (\sum_j \|S^N_j\|_{h^N}^2)^{1/2}$, so that
$|u_N|\equiv 1$.  For $c=(c_1,\dots,c_{d_N})$, we write
\begin{eqnarray*}\langle c,u_N(z)
\rangle &=& \left\langle c,\,\frac 1{|\Psi_N(z)|}\,\Psi_N
(z)\right\rangle\ =\
\frac 1{|\Psi_N(z)|}\,\sum_{j=1}^{d_N}c_j\,S^N_j(z)\in L_z^N\,,\\ |\langle
c,u_N(z)
\rangle|&=& \|\langle c,u_N(z)\rangle \|_{h^N}\,.\end{eqnarray*}

\begin{lem} \label{varinta}
$$\Var( Z_{s^N}) = -\frac{1}{\pi^2}
\d_z\dbar_z \d_w\dbar_w \int_{\C^{d_N}} \log |\langle c,u_N(z)\rangle|
\,\log |\langle c, u_N(w)\rangle|\, d\ga_N(c)\;.$$
\end{lem}

\begin{proof}

We write sections $s^N\in H^0(M,L^N)$ as
\begin{equation}\label{Psi0}s^N=\sum_{j=1}^{d_N} c_jS_j^N = \langle
c, \Psi_N\rangle\;,\qquad c=(c_1,\dots,c_{d_N}) \;.\end{equation}
Writing $\Psi_N=Fe_L^{\otimes N}$, where $e_L$ is a local
nonvanishing section of $L$, and recalling that $$\om=\frac
i2\Theta_h=- i \ddbar\log \|e_L\|_h\;,$$ we have by \eqref{Zs},
\begin{eqnarray}\label{ZsN}Z_{s^N}&=&\frac i\pi \ddbar \log
|\langle c,F\rangle| \ =\ \frac i\pi \ddbar \log |\langle
c,\Psi_N\rangle| -\frac i\pi \ddbar \log\|e_L^{\otimes
N}\|_h\nonumber\\& =& \frac i\pi \ddbar \log |\langle
c,\Psi_N\rangle|+ \frac N\pi \,\om\;.\end{eqnarray}

Consider the random current
\begin{equation}\label{YN}\wh Z_N:=\frac i\pi \ddbar \log
|\langle c,\Psi_N\rangle| = Z_{s^N} - \frac N\pi
\om\,.\end{equation} It follows immediately from the definition
\eqref{vc} of variance currents that
$$\Var (\wh Z_N) = \E(\wh Z_N\boxtimes \wh Z_N) -  \E(\wh Z_N)\boxtimes
\E(\wh Z_N) =\Var (Z_{s^N})\;.$$ By \eqref{Egamma}, we have
\begin{equation}\label{EYN} \E(\wh Z_N) =  \frac i\pi \ddbar \log
|\Psi_N| \;,\end{equation} whereas by \eqref{YN}, we have
\begin{eqnarray}\label{EZV} \E( \wh Z_N\boxtimes \wh Z_N)&=& -\frac{1}{\pi^2}
\int_{\C^{d_N}}\d_z\dbar_z \d_w\dbar_w   \log |\langle
c,\Psi_N(z)\rangle|\, \log |\langle c,\Psi_N(w) \rangle|\,
d\ga_N(c)\nonumber\\ &=& -\frac{1}{\pi^2}\, \d_z\dbar_z
\d_w\dbar_w \int_{\C^{d_N}}\log |\langle c,\Psi_N(z)\rangle|\,
\log |\langle c,\Psi_N(w) \rangle|\, d\ga_N(c)\;.
\end{eqnarray}
Recalling \eqref{uN}, we have
\begin{eqnarray}\log |\langle \Psi_N(z), c\rangle| \,\log |\langle
\Psi_N(w), c\rangle| &=&
\log |\Psi_N(z)| \,\log |\Psi_N(w)| + \log|\Psi_N(z)| \,\log |\langle c,
u_N(w)\rangle|\nonumber \\&&+ \log |\Psi_N(w)| \,\log |\langle c,
u_N(z)\rangle |\nonumber \nonumber\\&&+
 \log |\langle c, u_N(w)\rangle |  \,\log |\langle c, u_N(z)\rangle
|\;,\label{4terms}\end{eqnarray} which decomposes (\ref{EZV}) into
four terms. By (\ref{EYN}), the first term contributes $$-\frac
1{\pi^2}\, \ddbar\log|\Psi_N(z)|\wedge
\ddbar\log|\Psi_N(w)|=\E(\wh Z_N)\boxtimes\E(\wh Z_N)\,.$$ The
$c$-integral in the second term is independent of $w$ and hence
the second term vanishes when applying $\d_w\dbar_w$. The third
term likewise vanishes  when applying $\d_z\dbar_z$. Therefore,
the fourth term gives the variance current $\Var( Z_{s^N})$.
\end{proof}

To complete the proof of Theorem \ref{BIPOT}, we use the following
probability lemma, which gives the
$c$-integral of Lemma \ref{varinta}:

\begin{lem} \label{varintb} Let $(Y_1,Y_2)$ be joint complex Gaussian
random variables with  mean 0 and  $\E(|Y_1|^2)= \E(|Y_2|^2)=1$.  Then
$$\E\big(\log |Y_1|\, \log |Y_2|\big) = G\big(\left|\E(Y_1\overline
Y_2)\right|\big)\;,$$ where $$ G(t):=  \frac{\ga^2}{4}-\frac 1{4}
\int_0^{t^2} \frac{\log(1-s)}{s}\,ds\;, \quad 0\le t\le
1\qquad (\ga = \mbox{Euler's constant})\,.$$\end{lem}

\begin{proof} By replacing $Y_1$ with $e^{i\al}\,Y_1$, we can assume without
loss of generality that $\E(Y_1\overline Y_2)\ge 0$.   We can write
\begin{eqnarray*} Y_1 &=& \Xi_1\;,\\ Y_2 &=& (\cos\theta)\,\Xi_1 +
(\sin\theta)\,\Xi_2\,,
\end{eqnarray*} where $\Xi_1,\Xi_2$ are independent joint complex Gaussian
random variables with mean 0 and variance 1, and $\cos\theta=\E(Y_1\overline
Y_2)$.  Then
\begin{equation}\label{Elog} \E(\log |Y_1| \log |Y_2|) =
G(\cos\theta)\;,\end{equation} where \begin{equation}\label{G}
G(\cos\theta)=
\frac{1}{\pi^2} \int_{\C^2}  \log
|\Xi_{1}|\,\log \left|\Xi_1\cos\theta + \Xi_2\sin\theta
\right|\,e^{-(|\Xi_1|^2 + |\Xi_2|^2)}\, d\Xi_1
\,d\Xi_2\;.\end{equation}

The computation of $G(\cos\theta)$ was essentially given in
\cite[\S4.1]{BSZ1}.  We repeat this computation here  for
the readers' convenience:
Write
$\Xi_1 = r_1 e^{i \al},\ \Xi_2= r_2 e^{i(\al+\phi)}$, so that \eqref{G}
becomes
\begin{equation*}
G(\cos\theta)=\frac{2}{\pi}\int_0^{\infty} \int_0^{\infty}\int_0^{2\pi}
r_1 r_2
e^{-(r_1^2 + r_2^2)}  \log r_1 \log |r_1 \cos\theta+ r_2 e^{i
\phi}\sin\theta|\, d\phi\, dr_1\, dr_2\;.\end{equation*}  Evaluating the
inner integral by Jensen's formula, we obtain
\begin{equation*}  \int_0^{2\pi} \log|r_1 \cos \theta + r_2 \sin \theta e^{i
\phi}| \,d\phi = \left\{ \begin{array}{ll}2\pi \log (r_1 \cos \theta) &
\mbox{for}\;\;r_2 \sin \theta \leq r_1 \cos \theta \\ & \\ 2\pi\log (r_2
\sin \theta) & \mbox{for} \;\;r_2 \sin \theta\geq r_1 \cos \theta
\end{array}
\right. \end{equation*} Hence \begin{equation*} G(\cos\theta) = 4
\int_0^{\infty} \int_0^{\infty} r_1 r_2 e^{-(r_1^2 + r_2^2)}  \log r_1
\log
\max ( r_1 \cos \theta , r_2 \sin \theta)\,dr_1 \,dr_2. \end{equation*}
We make the change of variables  $r_1 = \rho \cos \phi,\, r_2 = \rho
\sin\phi$
to get \begin{equation*} G(\cos\theta) = 4 \int_0^{\infty}
\int_0^{\pi/2} \rho^3 e^{-\rho^2} \log (\rho \cos \phi) \log\max (\rho
\cos \phi \cos \theta , \rho \sin \phi \sin \theta)   \cos\phi\sin \phi
\,d\phi \,d\rho\;. \end{equation*} Since $$ \log \max( \rho \cos \phi \cos
\theta
, \rho \sin \phi \sin \theta) = \log (\rho \cos \phi \cos \theta) +
\log^+
(\tan \phi \tan \theta)\;,$$ we can write
$G=G_1+G_2$, where
\begin{eqnarray} \label{G1} G_1(\cos\theta)&=&4
\int_0^{\infty} \int_0^{\pi/2} \rho^3 e^{-\rho^2}   \log (\rho \cos
\phi)
\log ( \rho \cos \phi \cos \theta)  \cos\phi \sin \phi \,d\phi
\,d\rho\;,\\ \label{G2} G_2(\cos\theta) &=&4\int_0^{\infty} \int_{\pi/2 -
\theta}^{\pi/2}\rho^3 e^{-\rho^2}   \log (\rho \cos \phi) \log ( \tan
\phi
\tan \theta) \cos\phi\sin \phi \,d\phi \,d\rho\;. \end{eqnarray}

{From} (\ref{G1}),
$G_1(\cos\theta)=C_1+C_2\log\cos\theta$.  Substituting
$$\cos\theta = e^{-\la}\;,$$ we obtain
\begin{equation}\label{d0}G_1(e^{-\la})= C_1- C_2
\la\;.\end{equation} We now evaluate $G_2$. Since the integrand in
(\ref{G2}) vanishes when $\phi=\pi/2 -\theta$, we have
$$\frac{d }{ d\la}G_2(e^{-\la})= 4 \left(\frac{d }{ d\la}\log\tan\theta
\right) \int_0^{\infty} \int_{\pi/2 -
\theta}^{\pi/2}\rho^3 e^{-\rho^2} \log (\rho \cos \phi) \cos\phi\sin
\phi
\,d\phi \,d\rho\;.$$  Since  $$\frac{d }{ d\la}\log\tan\theta
=\frac 12\,\frac{d }{ d\la}\log(e^{2\la}-1)
=\frac{1}{
1-e^{-2\la}}\;,$$
we have $$\frac{d }{ d\la}G_2(e^{-\la})=  \frac{4}{
1-e^{-2\la}}\;(I_1+I_2)\;,$$ where
\begin{eqnarray*}I_1 &=&\int_0^{\infty} \int_{\pi/2 -
\theta}^{\pi/2}\rho^3 e^{-\rho^2} (\log \rho) \cos\phi\sin \phi
\,d\phi \,d\rho\ =\ C_3 \sin^2 \theta\ =\ C_3(1-e^{-2\la})\;,\\
I_2 &=&\int_0^{\infty} \int_{\pi/2 -
\theta}^{\pi/2}\rho^3 e^{-\rho^2} (\log \cos\phi) \cos\phi\sin \phi
\,d\phi \,d\rho\\ &=& \frac 12\int_{\pi/2 -
\theta}^{\pi/2} (\log \cos\phi) \cos\phi\sin \phi
\,d\phi \ =\ \frac 12\int_0^{\sin\theta}t\log t\,dt\\
&=&\frac{1}{ 8}(\sin^2\theta\log\sin^2\theta-\sin^2\theta)\ =\ \frac{1}{
8}(1-e^{-2\la})\left[\log (1-e^{-2\la})-1\right]\;. \end {eqnarray*}
Thus \begin{equation} \label{d1} \frac{d }{ d\la}G_2(e^{-\la})=
\frac{1}{
2}\log (1-e^{-2\la}) +4\,C_3 -\frac 12\;.\end{equation}

Combining \eqref{d0}--\eqref{d1}, we have
\begin{equation}\label{almostG} G(e^{-\la})=C_4+C_5\, \la +\frac 12
\int_0^\la \log(1-e^{-2s})\,ds\;.\end{equation} By \eqref{G},
$$G(0)=\big[\E(\log |\Xi_1|)\big]^2 = \left[2\int_0^\infty (\log
r)\,e^{-r^2}\,r\,dr\right]^2=
\frac{\ga^2}{4}\;.$$ Substituting $\la=\infty$ in
\eqref{almostG}, we conclude that $C_5=0$ and
\begin{equation}\label{finishG}
G(e^{-\la})=\frac{\ga^2}{4}-\frac 12 \int_\la^\infty \log(1-e^{-2s})\,ds\;,
\end{equation} or equivalently,
\begin{equation}\label{finishG1}
G(t)=\frac{\ga^2}{4} -\frac 14 \int_0^{t^2} \frac{\log(1-s)}{s}\,ds \qquad
(0\le t\le 1)\;.
\end{equation}
\end{proof}

\medskip\noindent{\it Proof of Theorem \ref{BIPOT}:\/} Fix points $z,w\in
M$, and let $x,y\in X$ with $\pi(x)=z,\ \pi(y)=w$. We apply Lemma
\ref{varintb} with $Y_1=\langle c, \hat u_N(x)\rangle$,  $Y_2=\langle c,
\hat u_N(y)\rangle$. Since $|\langle c, \hat
u_N\rangle|=  |\langle c,u_N\rangle|\circ\pi$, we have $$\log|Y_1| =  \log
|\langle c,u_N(z)\rangle|\,,\quad \log|Y_2| =  \log
|\langle c,u_N(w)\rangle|\,.$$
To determine $\E(Y_1\overline Y_2)$, we note that for a random $\hat s^N =
\sum c_j\wh S^N_j
\in \hcal^2_N(X)$,
\begin{equation} \label{EPi} \E\left(\hat s(x)\,\overline{\hat s(y)}\right)
= \sum_{j,k=1}^{d_N} \E(c_j\bar c_k) \,\wh S^N_j(x)\,\overline{\wh
S^N_k(y)} = \sum_{j=1}^{d_N}\wh S^N_j(x)\,\overline{\wh
S^N_j(y)} =\Pi_N(x,y)\;.\end{equation}
Since $$\langle c,\hat u_N(x)\rangle = \frac {\langle c,\wh\Psi_N(x)
\rangle} {\big|\wh \Psi_N(x)\big|}  = \frac {\hat
s^N(x)}{\Pi_N(x,x)^{1/2}}\,,$$ we have by \eqref{EPi},
$$\E(Y_1\overline Y_2) = \frac
{\Pi_N(x,y)}{\Pi_N(x,x)^{1/2}\Pi_N(y,y)^{1/2}}\,,$$ and
recalling \eqref{PN},
\begin{equation}\label{EYY} \big|\E(Y_1\overline Y_2)\big| =
P_N(z,w)\,.\end{equation}

Therefore, by Lemma \ref{varintb} and \eqref{EYY},
$$\int_{\C^{d_N}}\log |\langle u_N(z), c\rangle| \,
\log |\langle u_N(w), c\rangle|\, d\ga_N(c) = \E \big(\log |Y_1|\, \log
|Y_2|\big) = G\big(P_N(z,w)\big)\,.$$
By \eqref{Gtilde} and \eqref{finishG1},
\begin{equation}\label{GGtilde} \wt G(t)= \frac 1{\pi^2}
\left[G(t)-\frac{\ga^2}4\right]\;,\end{equation} and hence, recalling
that $\wt G\circ P_N=Q_N$,
\begin{equation}\label{equalsG}\frac 1{\pi^2}\int_{\C^{d_N}} \log |\langle
u_N(z), c\rangle| \,\log |\langle u_N(w), c\rangle|\, d\ga_N(c)=
Q_N(z,w)+C\;.\end{equation}
Theorem \ref{BIPOT} follows by combining
Lemma \ref{varinta} and \eqref{equalsG}.\qed

\subsection{Asymptotics of the pluri-bipotential}\label{Qasympt}

  We now use the \szego kernel off-diagonal asymptotics to describe the
$N$-asymptotics for the variance current $\Var(Z_{s^N})$ (Lemma
\ref{d4Qas}).  We also need to know the behavior of the variance
current near the diagonal.  We showed in \cite[(107)]{BSZ2} that
the codimension-one scaling limit pair correlation $K^\infty_{21}$
grows like\break $|z-w|^{-2}$ near the diagonal (for dimension
$m\ge 2$). Our computation of the variance current asymptotics
also gives this growth rate for the variance current (Lemma
\ref{DISTM2}) as well as for its scaling limit (Lemma
\ref{d4Qas}).

We begin by noting that the pluri-bipotential decays rapidly away from the
diagonal:

\begin{lem}\label{Qdecay}  For $b>\sqrt{j+q+1},\ j\ge 0$, we have
$$|\nabla^j Q_N( z,w)| =O\left(\frac
1{N^{q}} \right)\;,\quad \mbox{for }\ \dist(z,w) \ge \frac{b\sqrt{\log
N}}{\sqrtn}\;.$$ \end{lem}

\begin{proof}  We recall from \eqref{QN} that $Q_N=\wt G\circ P_N$,  where
$\wt G$ is analytic at 0 (with radius of convergence 1) and $\wt
G(t)=O(t^2)$.  The estimate then follows from
Proposition \ref{DPdecay} with $k = \left\lfloor
\frac{q+1}2\right\rfloor$.\end{proof}

Applying Lemma \ref{Qdecay} to the pluri-bipotential formula for
the variance of Theorem \ref{BIPOT}, we conclude that  the
variance current decays rapidly away from the diagonal.

We next show the near-diagonal estimate:
\begin{lem}\label{FNnear} For $b\in \R^+$, we have $$ Q_N\Big(z_0,z_0+\frac
v {\sqrtn}\Big)= \wt G(e^{-\frac 12 |v|^2})+
O(N^{-1/2+\epsilon})\;,\qquad \mbox {for }\ |v|\le b\sqrt{\log N
}.$$\end{lem}

\begin{proof}  Since
$P_N(z_0,z_0)=1$ and $\wt G'(t)
\to \infty$ as $t\to 1$, we need a short argument: let
\begin{equation}\label{Lambda} \Lambda_N= -\log
P_N\;.\end{equation} Recalling \eqref{Gtilde1}, we write,
\begin{equation}\label{F} F(\la):=\wt G(e^{-\la}) = -\frac
1{2\pi^2} \int_\la^\infty \log(1-e^{-2s})\,ds \qquad\quad
(\la\ge 0)\;,\end{equation} so that \begin{equation} \label{USEFUL}
Q_N=F\circ \Lambda_N\;. \end{equation}  By Proposition \ref{best},
\begin{equation}\label{RN0}\Lambda_N\left(z_0,z_0+\frac v{\sqrtn}\right)
={\half |v|^2} + \wt R_N(v) \;,\end{equation} where
\begin{equation}\label{RN}\wt
R_N=-\log(1+R_N)=O(|v|^2N^{-1/2+\ep})\quad \mbox{for }\
|v|<b\sqrt{\log N}\;.\end{equation} By \eqref{F},
\begin{equation}\label{F'} 0<-F'(\la) = -\frac
1{2\pi^2}\log(1-e^{-2\la}) \le \frac 1{2\pi^2}\left(1+\log^+
\frac1{\la}\right)\;.
\end{equation}
Since  $\half |v|^2+\wt R_N(v) = |v|^2\left(\half
+o(N)\right)$, it follows from
 \eqref{RN0}--\eqref{F'} that
\begin{eqnarray*}  Q_N\Big(z_0,z_0+\frac v
 {\sqrtn}\Big)&=&F\left(\half |v|^2+\wt R_N(v)\right)\nonumber\\&
= & F\left(\half |v|^2\right) +O\left(\left[1+\log^+
\frac1{|v|}\right]\wt R_N(v) \right)\nonumber \\&=& \wt
G(e^{-\frac 12 |v|^2})+ O(N^{-1/2+\epsilon})\;,\qquad \mbox {for
}\ |v|\le b\sqrt{\log N }\;.\end{eqnarray*}\end{proof}

We shall use the following notation: for a current $T$ on $M\times
M$, we  write
$$\d T=\d_1 T+\d_2 T\,,\qquad \d_1 = \sum dz_j\; \frac {\d}{\d z_j }\;,
\quad \d_2 = \sum dw_j\; \frac {\d}{\d w_j }\;,$$ where
$z_1,\dots,z_m$ are local coordinates on the first factor, and
$w_1,\dots, w_m$ are local coordinates on the second factor of
$M\times M$.  We similarly write
$$ \dbar T = \dbar_1 T +\dbar_2T\;.$$  In particular, we shall write
$\d_1\dbar_1\d_2\dbar_2Q_N$ in place of $\d_z\dbar_z \d_w\dbar_w
Q_{N}(z,w)$  to avoid confusion when we change variables.

 Next we
compute the leading term of the $N$-asymptotics of
$\dbar_1\dbar_2Q_N$ and $\d_1\dbar_1\d_2\dbar_2Q_N$. We choose
normal coordinates at a point $z_0\in M$, and we recall that in
terms of these coordinates, we have
\begin{equation}\label{Q}Q_N=F\circ\La_N,\quad
\La_N(z,w)=\frac N2|w-z|^2+\wt R_N(\sqrtn\,z,\sqrtn\,w)\;,\end{equation}
where
$\wt R_N$ is given by \eqref{RN}.
We now write $A_N(z,w)\approx B_N(z,w)$ when $$A_N(z,w)- B_N(z,w)=
O\left(N^{-\frac12+\ep}\,|B_N(z,w)|\right) \quad \mbox{for }\ |z|+|w| <
b\,\sqrt{\frac {\log N}{N}}\;.
$$ By \eqref{Q}, we have:   \begin{eqnarray}\dbar_2 Q_N(z,w)&\approx &\textstyle \frac N2
\,F'(\La_N(z,w))\,\dbar\big[(w-z)\cdot \bar w\big]\,,\label{d1Q}\\
\dbar_1\dbar_2 Q_N(z,w)&\approx &\textstyle -\frac {N^2}4
\,F''(\La_N(z,w))\, \dbar\big[(w-z)\cdot \bar z\big]\wedge \dbar
\big[(w-z)\cdot
\bar w\big]\,,\label{d2Q}\\ \dbar_1\d_2\dbar_2 Q_N(z,w)&\approx
&\textstyle -\frac {N^3}8 \,F^{(3)}(\La_N(z,w))\, \dbar\big[(w-z)\cdot
\bar z\big]\wedge \d \big[(\bar w-\bar z)\cdot w\big] \wedge
\dbar\big[(w-z)\cdot \bar w\big]\nonumber\\&&\textstyle\hspace{-1in} -  \frac
{N^2}4
\,F''(\La_N(z,w))\, \left\{\dbar\d\big[\bar z\cdot w\big]\wedge \dbar
\big[(w-z)\cdot \bar w\big] + \dbar\big[(w-z)\cdot \bar
z\big]\wedge \ddbar|w|^2\right\},\label{d3Q} \end{eqnarray}
and hence \begin{eqnarray}\label{d4Q}
\d_1\dbar_1\d_2\dbar_2Q_N(0,w) &\approx& \textstyle \frac
{N^4}{16}F^{(4)}(\La_N(0,w))\ \d (z\cdot\bar w)\wedge \dbar(\bar z\cdot w)
\wedge \d|w|^2\wedge \dbar|w|^2\nonumber \\&&\hspace{-1in}+
\textstyle \frac{N^3}{8} F^{(3)}(\La_N(0,w))\, \big[\ddbar|z|^2\wedge 
\d|w|^2\wedge \dbar|w|^2 +\dbar(\bar z\cdot w)\wedge \d|w|^2\wedge
\ddbar(z\cdot
\bar w)\nonumber\\&&\hspace{-1in}\textstyle
\qquad\qquad + \d(z\cdot \bar w)\wedge  \dbar \d (\bar z\cdot w)\wedge
\dbar|w|^2 + \d(z\cdot \bar
w)\wedge \dbar( \bar z\cdot w)\wedge \ddbar|w|^2
\big]\nonumber\\&&\hspace{-1in}+\textstyle \frac {N^2}4
\,F''(\La_N(0,w))\,\big[ \dbar \d (\bar z\cdot w)\wedge\ddbar(z\cdot\bar w) +
\ddbar|z|^2\wedge \ddbar |w|^2\big]\,.\end{eqnarray}  Differentiating \eqref{F'}, we have
\begin{equation}\label{F''}F''(\la)= \frac 1{\pi^2} \,\frac
1{e^{2\la}-1}\;,\end{equation} and hence
\begin{equation}\label{dF} F^{(j)}(\la)=O(e^{-2\la}) \qquad
(\la>1),
\end{equation}for $j\ge 0$.
Furthermore, by \eqref{F'},
$$F'(\la) = \frac 1{2\pi^2}\log\la + \eta(\la), \qquad
\eta\in\ccal^\infty([0,+\infty))\;,$$ and therefore
\begin{equation}\label{d4F}F^{(j+1)}(\la) = (-1)^{j+1} \frac
{(j-1)!}{2\pi^2}\, \la^{-j} +O(1) \qquad (\la>0),
\end{equation}for $j\ge 1$.

We now use the above computation  to describe the singularity of the
variance current near the diagonal.  We first recall an elementary fact:

\begin{lem}\label{elem} Let $u\in \ccal^1((\R^p\sm\{0\})\times
\R^q)$, and let $1\le j\le p+q$.  Suppose that $\d
u/\d x_j\in L^1(\R^{p+q})$ and
$u(x)=o\left(|\pi_1(x)|^{-p+1}\right)$, where $\pi_1:\R^{p+q}\to\R^p$ is
the projection.  Then the distribution derivative   $\d
u/\d x_j\in\dcal'(\R^{p+q})$ is given by the pointwise derivative, i.e.
$$\int_{(\R^p\sm\{0\})\times
\R^q} u\, \frac{\d\phi}{\d x_j} = -\int_{(\R^p\sm\{0\})\times
\R^q}
\frac{\d u}{\d x_j}\,\phi \qquad \forall\ \phi\in\dcal(\R^{p+q})\;.$$
\end{lem}
\begin{proof}Let $U_\ep=\{x\in\R^{p+q}:|\pi_1(x)|>\ep\}$.
The identity follows by integrating $u\, \frac{\d\phi}{\d x_j}$ by
parts over $U_\ep$, and noting that the boundary term
$$\int_{\d U_\ep} u\phi \,dx_1\cdots dx_{j-1}dx_{j+1}\cdots dx_{p+q}$$
goes to zero as $\ep\to 0$.\end{proof}

\begin{lem} \label{DISTM2} There   exist  a constant $C_m\in\R^+$
(depending only on the dimension $m$) and   an integer  $N_0=N_0(M)\in\Z^+$ such that for
$N\ge N_0$, we have:
\begin{itemize} \item[i)] The coefficients of the current
$\dbar_1\dbar_2 Q_N $ are locally bounded functions (given by
pointwise differentiation of $Q_N$), and we have the pointwise estimate
 $$|\dbar_1\dbar_2 Q_N(z,w) |\le C_mN\quad \mbox {for }\ 0<|w-z|<
b\sqrt{\frac {\log
 N}{N}}.
$$

\item[ii)] If $m\ge 2$, the  coefficients of the current
$\d_1\dbar_1\d_2\dbar_2Q_N$ are locally $L^{m-1}$ functions, and we
have the estimate
 $$|\d_1\dbar_1\d_2\dbar_2Q_N(z,w)| \le \frac{C_mN}{|w-z|^2}
\quad \mbox {for }\ 0<|w-z|< b\sqrt{\frac {\log N}{N}}.
$$
\end{itemize}\end{lem}

\begin{proof}
We take $z=z_0=0$.  By Propositions \ref{better}--\ref{best}, we
can choose $N_0$ such that
\begin{equation}\label{lowerlambda}\La_N(0,w) \ge \frac N3 |w|^2
\quad\mbox{for }\ |w|< b\,\sqrt{\frac {\log N}{N}}\ , \quad N\ge
N_0\;.\end{equation} By applying the chain rule as in
\eqref{d1Q}--\eqref{d4Q}, we conclude that for each $N\ge N_0$,
\begin{eqnarray}\nabla Q_N(z,w) \ =\ O(|w-z|\log |w-z|),\ &&\ \nabla^2
Q_N(z,w)\ =\ O(\log|w-z|),\nonumber\\
\nabla^j Q_N(z,w) &=& O(|w-z|^{-j+2}) \quad \mbox{for }\ j\ge
3.\label{ddd}\end{eqnarray} Hence, the partial derivatives of $Q_N$ of order
$\le 3$ are in $L^1_{loc}$, and the same holds for the fourth
order derivatives if $m\ge 2$.  By repeatedly applying Lemma
\ref{elem} with $\pi_1(x)=w-z$, we conclude that the currents
$\dbar_1\dbar_2 Q_N $ and $\d_1\dbar_1\d_2\dbar_2Q_N$ have locally
$L^1$ coefficients.  The upper bound in (i) follows from \eqref{d2Q},
\eqref{F''} and \eqref{lowerlambda}, and hence the coefficients of
$\dbar_1\dbar_2 Q_N
$ are actually in $L^\infty_{loc}$. The upper bound in (ii) similarly
follows from (\ref{d4Q}),
(\ref{d4F}) and \eqref{lowerlambda}, and hence the coefficients of
$\d_1\dbar_1\d_2\dbar_2Q_N$ are in $L^{m-1}_{loc}$.
\end{proof}

\medskip These computations show that $Q_N$ is $\ccal^1$ and has
vanishing first derivatives on the diagonal in $M\times M$, but
$Q_N$ is not $\ccal^2$ along the diagonal.  Lemma \ref{DISTM2}(i)
says that $\dbar_1\dbar_2Q_N$ is bounded; however, a similar
computation shows that $\d_2\dbar_2Q_N(z,w) \ge c \log|z-w|$, for
a positive constant $c$.  When $m>1$, $\d_1\dbar_1\d_2\dbar_2Q_N
\sim |w-z|^{-2}$; but when $m=1$, $\d_1\dbar_1\d_2\dbar_2Q_N$ is a
measure with a singular component along the diagonal, and off the
diagonal there is cancellation in \eqref{d4Q} and
$\d_1\dbar_1\d_2\dbar_2Q_N \sim |w-z|^{2}$ (see
\cite[Th.~4.2]{BSZ1}).

Making the change of variables
$$w = \frac v\sqrtn$$ in \eqref{d2Q} and again applying \eqref{Q} and
\eqref{d4F}, we obtain the following asymptotic formulas:

\begin{lem}\label{d2Qas}For $N$ sufficiently large,   
$$\textstyle\dbar_1\dbar_2 Q_N(z_0, z_0+\frac
v{\sqrtn}) = -\frac{\sqrtn}{4} F''(\half |v|^2)\,\dbar(\bar z\cdot v)
\wedge \dbar|v|^2 +O(N^\ep)\;,$$ for $0<|v|<b\sqrt{\log N}$. \end{lem}
\begin{proof}  By \eqref{RN0}--\eqref{RN} and \eqref{d4F} with $j=2$, we
have
\begin{eqnarray}\textstyle F''\left(\La_N(0,\frac v{\sqrtn})\right)&=&
F''(\half |v|^2)+ F^{(3)}\left(\left[\half
+O(N^{-1/2+\ep})\right]|v|^2\right)\cdot
O\left(|v|^2N^{-1/2+\ep}\right)\nonumber\\ &=& F''(\half |v|^2)+
O\left(|v|^{-2}N^{-1/2+\ep}\right).\label{mvt}
\end{eqnarray} The formula follows from \eqref{d2Q} and \eqref{mvt}.
\end{proof}

\begin{lem}\label{d4Qas}For $N$ sufficiently large,
\begin{equation}\label{d4Qv} \textstyle-
\d_1\dbar_1\d_2\dbar_2Q_N(z_0,z_0+ \frac v{\sqrtn}) =
N\,\Var_\infty^{z_0}(v) +O\left(|v|^{-2}\,N^{1/2+\ep}\right)
\quad\mbox{for }\ 0<|v|<b\sqrt{\log N} \;,\end{equation} where
$\Var_\infty^{z_0}\in   T^{*2,2}_{(z_0,v)}(M\times \C^m)$ is given by   
\begin{eqnarray}\label{many} \textstyle\Var_\infty^{z_0}(v) &:=&
\textstyle -\frac
{1}{16}\,F^{(4)}(\half |v|^2)\ \d (z\cdot\bar v)\wedge \dbar(\bar z\cdot
v)
\wedge \d|v|^2\wedge \dbar|v|^2\nonumber \\&&
\textstyle -\frac{1}{8}\, F^{(3)}(\half |v|^2)\, \big[\ddbar|z|^2\wedge 
\d|v|^2\wedge \dbar|v|^2 +\dbar(\bar z\cdot v)\wedge \d|v|^2\wedge
\ddbar(z\cdot
\bar v)\nonumber\\&&\textstyle
\qquad\qquad + \d(z\cdot \bar v)\wedge  \dbar \d (\bar z\cdot v)\wedge
\dbar|v|^2 + \d(z\cdot \bar
v)\wedge \dbar( \bar z\cdot v)\wedge \ddbar|v|^2
\big]\nonumber\\&& \textstyle -\frac {1}4
\,F''(\half |v|^2)\,\big[ \dbar \d (\bar z\cdot v)\wedge\ddbar(z\cdot\bar v) +
\ddbar|z|^2\wedge \ddbar |v|^2\big]\,.\end{eqnarray} Furthermore,
\begin{equation}\label{varest}
\Var_\infty^{z_0}(v)=\left\{\begin{array}{ll} O(|v|^{-2})\quad
&\mbox{for }\ |v|>0\\O(|v|^4\,e^{-|v|^2}) &\mbox{for }\
|v|>1\end{array}\right.\ .\end{equation}
\end{lem}
\begin{proof} Formula \eqref{many}
follows by the same argument as in the proof of Lemma \ref{d2Qas},
applying \eqref{d4Q} in place of \eqref{d2Q}. The estimate
\eqref{varest} follows by applying \eqref{dF}--\eqref{d4F} to
\eqref{many}.
\end{proof}

The current $\Var_\infty^{z_0}$ has $L^1_{loc}$ coefficients if
$m\ge 2$, but  contains the singular term   $\frac 1\pi
\om(z_0)\otimes \de_{z_0}(v)$
 if
$m=1$ (see   Theorem \ref{variant1} or \cite[Theorem~4.1]{BSZ1}). 

\subsection{The pair correlation current}
\label{highcodim}

The pair correlation current gives the correlation for the zero
densities at two points of $M$. It is defined to be
\begin{equation}\label{pcc} \K_{2k}^N:=
\E\big(Z_{s_1^N,\dots,s_k^N} \boxtimes Z_{s_1^N,\dots,s_k^N}
\big)\in \dcal'^{2k,2k}(M\times M)\;.\end{equation}  Thus by Corollary
\ref{indep},
\begin{equation}\label{vark}\Var\big(Z_{s_1^N,\dots,s_k^N}\big)=\K_{2k}^N-
\E\big(Z_{s_1^N,\dots,s_k^N}\big)\boxtimes\E
\big(Z_{s_1^N,\dots,s_k^N}\big) = \K_{2k}^N-
\left[\E\big(Z_{s^N}\big)\boxtimes\E\big(Z_{s^N}\big)\right]^k\; .\end{equation}
As a consequence of Theorem \ref{BIPOT}, we have the following
formula for the case $k=1$:
\begin{prop} \label{KN21} The pair correlation
current in codimension one is given by
$$\K^N_{21} = - \d_1\dbar_1 \d_2\dbar_2 Q_{N} +\E Z_{s^N}\boxtimes
\E Z_{s^N}\;.$$ For $m\ge 2$, the coefficients of the current $\K^N_{21}$ are in
$L^{m-1}_{loc}$.\end{prop}

\begin{proof} The formula for $\K^N_{21}$ is an immediate
consequence of Theorem \ref{BIPOT} and \eqref{vark}. By Lemma
\ref{DISTM2}(ii), the coefficients of the current
$\d_1\dbar_1\d_2\dbar_2Q_N$ are in $L^{m-1}_{loc}$, and hence the
same holds for the coefficients of $\K^N_{21}$.\end{proof}

We note that the pair correlation function $K_{2k}^N(z,w)$, which
gives the probability density of zeros occurring at both $z$ and
$w$ (see \cite{BSZ1, BSZ2}), can be obtained from the pair
correlation current: \begin{equation}\label{pcf}\textstyle
K_{2k}^N \left(\frac 1{m!}\,\om^m \boxtimes \frac
1{m!}\,\om^m\right)= \K_{2k}^N \wedge \left(\frac
1{(m-k)!}\,\om^{m-k} \boxtimes \frac
1{(m-k)!}\,\om^{m-k}\right)\;.\end{equation}  The advantage of the
pair correlation current is that, because of the independence of
the $s^N_j$, the codimension-$k$ correlation current is  the $k$-th exterior
power of the corresponding codimension-one current, i.e.
\begin{equation}\label{powers}\K_{2k}^N = \left[\K_{21}^N
\right]^{\wedge k}.\end{equation}
This formula is analogous to the corresponding identity
$\E(Z_{s^N_1,\dots ,s_k^N})=\E(Z_{s^N})^k$  for the expected value in
Corollary \ref{indep}; both formulas hold since the
$s^N_j$ are independent random sections.  However, in the case of
correlation currents, the right side of \eqref{powers} is not well defined
along the diagonal in $M\times M$,
since $\K^N_{21}$ is singular on the diagonal. We shall
show that for
$k<m$, the current
$\K^N_{2k}$ has
$L^1$ coefficients and \eqref{powers} holds, with $\left[\K_{21}^N
\right]^{\wedge k}$ given by pointwise multiplication.   However, for the
point case $k=m$ of Theorem \ref{number}, the pair
correlation current $\K^N_{2m}$ contains a singular measure supported on the
diagonal (see Theorem \ref{variant1}), and the right side of \eqref{powers} must be
interpreted as a limit of smooth currents.  The singularities of
$\K_{21}^N$ necessitate a more complicated proof, using a
smoothing method.  (One can also define in an analogous way the
$n$-point correlation currents $\K_{nk}^N$, which satisfy the identity
$\K_{nk}^N = \left[\K_{n1}^N\right]^{\wedge k}$.)

\subsection{Explicit formula for the variance}\label{explicit}

We shall write
\begin{equation}\label{Phi}\begin{array}{rcl}\Phi_k &=&
\frac 1 {(m-k)!}\om^{m-k} \quad \mbox{for }\ 1\le k\le
m-1,\\[8pt] \Phi_m &=&1\;,\end{array}\end{equation} so that
$$\vol_{2m-2k}\left(Z_{s_1^N,\dots,s_k^N}\cap U\right) =
\left(Z_{s_1^N,\dots,s_k^N}\,,\chi_U\,\Phi_k\right)\;, \quad \mbox{for }\
1\le k\le m\;.$$
By Corollary \ref{indep}, the expected volume of the codimension-$k$ zero current is
given by
\begin{equation} \E\big(\vol_{2m-2k}[Z_{s_1^N,\dots,s_k^N}\cap U]\big) =
\int_U \big(\E Z_{s^N}\big)^k \wedge \Phi_k\;,\end{equation} where
\begin{equation}\label{EZN} \E Z_{s^N}=   \frac i{2\pi} \ddbar
\log \Pi_{N}(z,z) + \frac N\pi \omega\;.\end{equation}

In this section, we prove the following   integral formula
for the volume and number variance:

\begin{theo}\label{varint2}The variance in Theorem \ref{volume} is given
by:
\begin{multline*}
\var\big(\vol_{2m-2k}[Z_{s_1^N,\dots,s_k^N}\cap U]\big)
\\=\sum_{j=1}^k(-1)^j {k \choose j}\int_{\d U\times
\d U}\dbar_1\dbar_2Q_N\wedge
\left(\d_1\dbar_1 \d_2\dbar_2
 Q_{N}\right)^{j-1}\wedge
\big(\E Z_{s^N} \boxtimes \E Z_{s^N}\big)^{k-j}\wedge
(\Phi_k\boxtimes\Phi_k),
\end{multline*} for $N$ sufficiently large, where $Q_N$ is given
by \eqref{QN}, $\E Z_{s^N}$ is given by \eqref{EZN} and the integrands  are
in
$L^1(\d U\times \d U)$. \end{theo}

In particular, for the one-dimensional case $k=m=1$, we have
$$\var\big(\ncal^U_N\big) =- \int_{\d U\times
\d U}\dbar_z\dbar_wQ_N(z,w)\;.$$

Theorem \ref{varint2} follows formally from Theorem
\ref{BIPOT} and equations \eqref{vark} and \eqref{powers}. To
verify the formula rigorously, we must show that the currents
$\dbar_1\dbar_2Q_N\wedge\big[\d_1\dbar_1 \d_2\dbar_2
Q_{N}\big]^{j-1}$ ($1\le j\le m$), which are smooth forms off the
diagonal in $M\times M$, are well defined and have $L^1$
coefficients; i.e., they impart no mass to the diagonal. To do this, we use
the asymptotics of $Q_N(z,w)$ as $|z-w|\to 0$ given in Lemma
\ref{DISTM2}.   In \S\ref{s-number}, we shall
use the $N$ asymptotics of Lemmas \ref{d2Qas}--\ref{d4Qas}
together with Theorem \ref{varint2} to prove
Theorem \ref{volume}.

\begin{defin} We say that a
current $u\in\dcal'^{p,q}(X)$ on a \kahler manifold $X$ is an $L^1$ current on $X$ if
its local coefficients are $L^1$ functions and $\int_X|u|\,d\vol_X<+\infty$.
(The second condition is redundant if $X$ is compact.)  If $\{u_n\}$ is
a sequence of $L^1$ currents on $X$, we say that $u_n\to u$ in $L^1$ if
$\int_X|u_n-u|\,d\vol_X \to 0$.\end{defin}

We shall prove  Theorem \ref{varint2} by approximating $\chi_U$ by smooth
cut-off functions, and then applying the following explicit formula for the
variance current:

\begin{theo}\label{variant} Let  $1\le k\le m$.
Then for $N$ sufficiently large,
\begin{multline*} \Var\big(Z_{s_1^N,\dots,s_k^N}\big)
= \d_1\d_2\left[\sum_{j=1}^k   (-1)^{j-1}  {k\choose j}\dbar_1\dbar_2
Q_N\wedge\left(\d_1\dbar_1 \d_2\dbar_2 Q_{N}\right)^{j-1}\wedge \big(\E
Z_{s^N} \boxtimes \E Z_{s^N}\big)^{k-j}\right].
\end{multline*}
where  the current inside the brackets  is an $L^1$ current on $M\times M$ given by
pointwise multiplication,
$Q_N$ is given by \eqref{QN}, and  $\E Z_{s^N}$ is given by \eqref{EZN}.
Furthermore,
$\Var\big(Z_{s_1^N,\dots,s_k^N}\big)$ is an $L^1$ current on $M\times M$ if $k\le
m-1$.
\end{theo}

We need the following smoothing result for our proof of Theorem
\ref{variant}:

\begin{prop} \label{converge} Let $(L,h)\to M$ be as in Theorem
\ref{number}, with $\dim M = m\ge 2$.  Then there is a positive integer $N_0$ such
that the following holds for all
$N \ge N_0$:  \begin{itemize} \item for
$1\le j\le m-1$, \
$\K^N_{2j}$ is an $L^1$ current on $M\times M$, and is
given by the pointwise formula $\K^N_{2j}=(\K^N_{21})^{j}$.
\item for all points $P_0\in M\times M$, there exist a neighborhood $\Om\subset
M\times M$ of $P_0$ and smooth forms
$$S_\ep K^N_{21}\in \ecal^{2,2}(\Om) \qquad(0<\ep<1)$$ such that
\begin{itemize}
\item[i) ] $\d_1(S_\ep K^N_{21})= \d_2(S_\ep K^N_{21}) = 0$;
\item[ii) ] for $1\le
j\le m-1$, \
$(S_\ep
\K^N_{21})^{j}
\to \K^N_{2j}|_\Om$ in $L^1$, as $\ep\to 0$;
\item[iii) ] for $2\le k\le m$,  $\K^N_{21} \wedge \left(S_\ep
\K^N_{21}\right)^{k-1}\to \K^N_{2k}|_\Om $ weakly, as
$\ep\to 0$.
\end{itemize}
\end{itemize}\end{prop}

We postpone the proof of Proposition \ref{converge} to the next section,
and we now prove Theorem \ref{variant}, assuming Proposition
\ref{converge}:  The case $m=1$ of Theorem \ref{variant} is Theorem \ref{BIPOT}.  So
we let $m\ge 2$. It suffices to consider a test form
$\phi\in\dcal^{2m-2k,2m-2k}(\Om)$, where $\Om\subset
M\times M$ is as in Proposition
\ref{converge}. By Propositions \ref{KN21} and
\ref{converge} (recalling that
$\dbar_1\dbar_2Q_N$ has $ L^\infty$ coefficients by Lemma \ref{DISTM2}), we have
\begin{eqnarray}(\K^N_{2k},\phi) &=& \lim_{\ep\to 0}
\big( \K^N_{21}\wedge [S_\ep \K^N_{21}]^{k-1},\,
\phi\big)\nonumber\\ & =& \lim_{\ep\to
0}\int_{\Om}\Big\{ \dbar_1\dbar_2Q_{N}\wedge
[S_\ep \K^N_{21}]^{k-1}\wedge\d_1\d_2\phi\ +\ [\E Z_{s^N}\boxtimes \E
Z_{s^N}]\wedge [S_\ep
\K^N_{21}]^{k-1}
\wedge\phi\Big\}\nonumber
\\ &=&
\int_{\Om}
 \dbar_1\dbar_2Q_{N} \wedge ( - \d_1\dbar_1 \d_2\dbar_2 Q_{N} +\E
Z_{s^N}\boxtimes \E Z_{s^N})^{k-1}\wedge\d_1\d_2
\phi\nonumber\\&& +\int_{\Om}(\E Z_{s^N}\boxtimes \E Z_{s^N})\wedge
( - \d_1\dbar_1 \d_2\dbar_2 Q_{N} +\E Z_{s^N}\boxtimes \E
Z_{s^N})^{k-1}\wedge\phi\;.\label{v1}
\end{eqnarray}
Expanding the integrand and recalling \eqref{vark}, we then
have
\begin{eqnarray*}\left(\Var\big(Z_{s_1^N,\dots,s_k^N}\big),\,\phi
\right)\hspace{-1in}\\  &=&\left(\K_{2k}^N,\,\phi\right)-
\int
(\E Z_{s^N}\boxtimes \E Z_{s^N})^{k}\wedge\phi\\
&=& \sum_{j=1}^{k} {k-1\choose j-1}\int \dbar_1\dbar_2 Q_{N} \wedge ( - \d_1\dbar_1
\d_2\dbar_2 Q_{N})^{j-1}\wedge (\E Z_{s^N}\boxtimes \E
Z_{s^N})^{k-j}\wedge\d_1\d_2\phi\\&& +
\sum_{j=1}^{k-1} {k-1\choose j}\int \dbar_1\dbar_2 Q_{N}\wedge  ( - \d_1\dbar_1
\d_2\dbar_2 Q_{N})^{j-1}\wedge (\E Z_{s^N}\boxtimes \E
Z_{s^N})^{k-j}\wedge\d_1\d_2\phi \\&=&
\sum_{j=1}^{k} {k\choose j}\int \dbar_1\dbar_2 Q_{N} \wedge  ( - \d_1\dbar_1
\d_2\dbar_2 Q_{N})^{j-1}\wedge (\E Z_{s^N}\boxtimes \E
Z_{s^N})^{k-j}\wedge\d_1\d_2\phi\,.\end{eqnarray*}
This is the formula of Theorem \ref{variant}; to complete the proof of
the theorem, it remains only to prove   Proposition
\ref{converge}.

\subsubsection{Off-diagonal decay of the correlation current}
We now use \ref{variant}--\ref{converge} to give a more explicit
formula for the pair correlation current in higher codimension and
describe its off-diagonal asymptotics.  (The results of this
section are presented here for their general interest and are not
needed for the proof of the variance formulas.)

We note that the correlation currents $\K^N_{2k}$ are smooth forms
away from the diagonal in $M\times M$.  We now show that
$\K^N_{2k}$ has no mass on the diagonal for $k<m$, while
$\K^N_{2m}$ contains a `delta-function' along the diagonal:

\begin{theo}\label{variant1}Let $(L,h)\to M$ be as in Theorem
\ref{number}.  Then  for
$N$ sufficiently large, we have:  \begin{itemize} \item[i)] for  $1\le k\le m-1$,
the correlation current for $\K^N_{2k}$ is an $L^1$ current on $M\times M$ given by
$$ \K^N_{2k}
= \left[ - \d_1\dbar_1 \d_2\dbar_2 Q_{N} +\E Z_{s^N}\boxtimes
\E Z_{s^N}\right]^k\;,
$$ where  $Q_N$
is given by \eqref{QN}, and  $\E Z_{s^N}$ is given by \eqref{EZN};\\[-8pt]
\item[ii)] $ \K^N_{2m}
=  (- \d_1\dbar_1 \d_2\dbar_2 Q_{N} +\E Z_{s^N}\boxtimes
\E Z_{s^N})^m\big|_{M\times M\sm\De} + \operatorname{diag}_* (\E
Z_{s^N})^m$\,,\\[4pt] where  $\De=
\{(z,z):z\in M\}$ is the diagonal in $M\times M$, and $\operatorname{diag}:M\to
M\times M$ is the diagonal map
$\operatorname{diag}(z)=(z,z)$.
\end{itemize} \end{theo}

\begin{proof} Part (i) is an immediate consequence of Proposition \ref{KN21} and the
first conclusion of Proposition \ref{converge}. To verify (ii),
we regard the current $\K^N_{2m}\in\dcal'^{4m}(M\times M)$ (which
is of order 0 by its definition) as a measure on $M\times M$.
Since $Q_N$ is $\ccal^\infty$ in $M\times M\sm\De$,  it follows
from Theorem \ref{variant} (or from \eqref{v1}) that $\K^N_{2m}
=(- \d_1\dbar_1 \d_2\dbar_2 Q_{N} +\E Z_{s^N}\boxtimes \E
Z_{s^N})^m$ on $M\times M\sm\De$.  Hence it suffices to show that
(ii) holds on $\De$, i.e., for any Borel set $A\subset M$,
\begin{equation} \label{diag} \big(\K^N_{2m}\,,\,
\chi_{\operatorname{diag}(A)}\big) =
\big((\E Z_{s^N})^m\,,\,\chi_A\big)\;.\end{equation} To verify \eqref{diag}, we note
that
$$\left(Z_{s^N_1,\dots,s^N_m}\boxtimes
Z_{s^N_1,\dots,s^N_m}\,,\,\chi_{\operatorname{diag}(A)} \right)=
\#\left\{z\in A: s^N_1(z)=\cdots=s^N_m(z)=0\right\}  =  \left(
Z_{s^N_1,\dots,s^N_m}\,,\,
\chi_A\right).$$ Taking  expectations and recalling Corollary \ref{indep}, we then
have
\begin{eqnarray*} \big(\K^N_{2m}\,,\, \chi_{\operatorname{diag}(A)}\big)
&=& \E \left(Z_{s^N_1,\dots,s^N_m}\boxtimes
Z_{s^N_1,\dots,s^N_m}\,,\,\chi_{\operatorname{diag}(A)} \right)\\&=& \left( \E
Z_{s^N_1,\dots,s^N_m}\,,\,
\chi_A\right)\ =\ \big((\E Z_{s^N})^m\,,\,\chi_A\big)\;,\end{eqnarray*}
and therefore (ii) holds on $\De$ and hence on all of $M\times M$.
\end{proof}

\begin{theo} \label{COR1decay} Let  $(L,h)\to M$ be as in
Theorem
\ref{number} and let $1\le k\le m$.
 For  $b>\sqrt{q+2k+3}\ge 3$,    the  pair correlation
currents satisfy the off-diagonal asymptotics
$$\K^N_{2k}   = \big(\E Z_{s^N}\boxtimes
\E Z_{s^N}\big)^k   + O\left(N^{-q} \right)\;,\quad \mbox{for }\
\dist(z,w)\ge\frac{b\sqrt{\log N}}{\sqrtn}\;.$$ In particular,
$$\K^N_{2k}   = \big(\E Z_{s^N}\boxtimes
\E Z_{s^N}\big)^k  + O\left({N^{-\infty}} \right)\;,\quad \mbox{for }\
\dist(z,w)\ge N^{-1/2+\ep}\;.$$
\end{theo}

\begin{proof} By  Lemma \ref{Qdecay}, for $b>\sqrt{q+2k+3}$\,, we
have
\begin{equation}\label{d4Qdecay} \d_1\dbar_1 \d_2\dbar_2
Q_{N} = O\left(\frac 1{N^{q+2k-2}} \right)\;,\quad \mbox{for }\ \dist(z,w)\ge
\frac{b\sqrt{\log N}}{\sqrtn}\;. \end{equation}
Since the statement only pertains to the
off-diagonal, by Theorem \ref{variant1} and \eqref{d4Qdecay},
\begin{equation}\label{byv}\K_{2k}^N = \left[ - \d_1\dbar_1 \d_2\dbar_2 Q_{N} +\E
Z_{s^N}\boxtimes
\E Z_{s^N}\right]^k = \left[ \E Z_{s^N}\boxtimes
\E Z_{s^N} + O\left(\frac 1{N^{q+2k-2}} \right)\right]^k,\end{equation} for
$\dist(z,w)\ge
\frac{b\sqrt{\log N}}{\sqrtn}$.  By
Proposition
\ref{indepTYZ}, $\E Z_{s^N}\boxtimes
\E Z_{s^N} = O(N^2)$, and the desired asymptotics
follow immediately from \eqref{byv}.\end{proof}

In particular, the point case of Theorem \ref{COR1decay} yields the asymptotics

\begin{cor}\label{cordecay}  Let  $(L,h)\to M$ be as in
Theorem
\ref{number}, and  define the normalized
pair correlation function by $$\wt K_{2m}^N  =
\frac{\K_{2m}^N}{\left(\E Z_{s^N}\boxtimes \E Z_{s^N}\right)^m}\;.
$$
Then for  $b>\sqrt{q+5}$, $q\ge 1$, we have
$$\wt K_{2m}^N(z,w)  = 1 + O\left(N^{-q}\right),\quad \mbox{for }\
\dist(z,w)\ge\frac{b\sqrt{\log N}}{\sqrtn}\;.$$ In particular,
$$\wt K_{2m}^N(z,w)  = 1 + O\left({N^{-\infty}} \right)\;,\quad \mbox{for }\
\dist(z,w)\ge N^{-1/2+\ep}\;.$$\end{cor}

Corollary \ref{cordecay} can also be obtained from
Theorem \ref{near-far} (or \cite[(95)--(96)]{SZ2}), using the argument in Section 4.1
of \cite{BSZ2}.

\subsection{Smoothing $\K^N_{21}$: Proof of Proposition
\ref{converge}}\label{s-smoothing}  We shall use the  following
fact about averaging currents of integration over a smooth family
$\{Y_t\}$ of submanifolds:

\begin{lem}\label{intersectEZ}Let $X$ and $ \Om$ be complex manifolds of
dimension $m$ and $n$ respectively, and let $X'$ be a complex submanifold
of
$X$. Let $\ycal$ be a  complex submanifold of $M\times
\Om$ such that the projections $\pi_1:\ycal\to M$ and $\pi_2:\ycal\to\Om$
are submersions, and let $Y_t=\pi_1\left(\pi_2\inv\{t\}\right) \subset X$,
for $t\in\Om$. Then for all
$\al\in\dcal^{n,n}(\Om)$, we have
\begin{itemize}
\item   $\displaystyle \int_{t\in\Om} [Y_t]\,\al(t)\in \ecal^{p,p}(X)$,
where $p=\codim\ycal = \codim_XY_t$\;;

\item for almost all $t\in\Om$, $X'\cap Y_t$ is a complex submanifold of
$X$ of codimension $p':=p+\codim X'$, and $$\int_{t\in\Om}[X'\cap Y_t]\,
\al(t) = [X']\wedge  \int_{t\in\Om} [Y_t^j]\,\al(t) \in
\dcal'^{p',p'}(X)\;.$$
\end{itemize} \end{lem}

\begin{proof} Let $\ycal'=\pi_1\inv(X')=\ycal\cap(X'\times\Om)$, and
consider the commutative diagram
\begin{equation}\label{diagram}\begin{array}{llllllll} & \ycal' & \buildrel
{\hat\iota}\over
\hookrightarrow &\ycal \\[6pt]
\pi_1'\swarrow & &  \pi_1\swarrow &&\searrow\pi_2 \qquad ,\\[6pt]
X' &  \buildrel {\iota}\over\hookrightarrow & X && \quad\ \Om
\end{array}\end{equation}
where $\pi_1'=\pi|_{\ycal'}$, and $\iota,\hat\iota$ are the inclusions.
Since $\pi_2$ is a submersion,
$Y_t$ is a smooth submanifold of $X$ for all $t\in\Om$.  For a test form
$\phi\in\dcal^{m-p,m-p}(X)$, we have
$$ \left(\int_{t\in\Om} [Y_t]\,\al(t)\,,\,\phi\right)\ \eqd\int_\Om
\al(t) \int _{Y_t} \phi = \int_\ycal \pi^*_2\al \wedge
\pi^*_1\phi = \left(\pi_{1*}\pi^*_2\al,\phi\right)\,,$$  and thus
\begin{equation}\label{pushforward}\int_{t\in\Om} [Y_t]\,\al(t)=
\pi_{1*}\pi^*_2\al\;,\end{equation} which is a smooth form, since the
push-forward of a smooth current via a submersion is smooth.

Since $\pi_1$ is a submersion, $\ycal'$ is a submanifold, and hence by
Sard's theorem, the set of critical values of the map
$\pi_2':=\pi_2|_{\ycal'}$ has measure zero.  Therefore $\pi'^{-1}_2(t) =
(X'\cap Y_t)\times \{t\} $ is a submanifold for almost all $t\in\Om$.

Now suppose that $\phi\in\dcal^{m-p',m-p'}(X)$.   Let $S$ denote the
set of critical points of the map
$\pi'_2:\ycal'\to \Om$, and let
$E=\pi'_2(S\cap\supp \pi'^*_1\phi)$, which is a closed subset of
measure zero in $\Om$. Let $\rho_n\in\ccal^\infty(\Om)$ such that
$\rho_n\ge 0$, $\rho_n\nearrow \chi_{\Om\sm E}$, and $\supp\rho_n\cap
E=\emptyset$. As before, we have
\begin{equation*}\int_\Om \rho_n(t)\al(t)
\int _{X'\cap Y_t} \phi =
\int_{\ycal'\sm S} \pi'^*_2(\rho_n\al) \wedge \pi'^*_1\phi =
\int_{\ycal'} \pi'^*_2(\rho_n\al) \wedge \pi'^*_1\phi\;.\end{equation*}
We claim that
\begin{eqnarray}\label{limit1}
\int_\Om \rho_n(t)\al(t)
\int _{X'\cap Y_t} \phi  &\to &
\int_\Om \al(t)
\int _{X'\cap Y_t} \phi \,,\\ \label{limit2}
\int_{\ycal'} \pi'^*_2(\rho_n\al) \wedge \pi'^*_1\phi &\to &
\int_{\ycal'} \pi'^*_2\al \wedge \pi'^*_1\phi\, ,\end{eqnarray}
as $n\to\infty$, and hence
\begin{equation}\label{limit3} \int_\Om \al(t)
\int _{X'\cap Y_t} \phi = \int_{\ycal'} \pi'^*_2\al \wedge
\pi'^*_1\phi\,.\end{equation} To verify \eqref{limit1}--\eqref{limit3},
 we first consider the case  $\al\ge 0$ and $\phi=\beta^{m-p'}$, where
$\beta$ is a (compactly supported) semi-positive (1,1)-form  on
$X$, and apply monotone convergence to obtain
\eqref{limit1}--\eqref{limit2} for this case.  Since the right side of
\eqref{limit3} is given by the integral of a compactly  supported smooth
form
 over $\ycal'$, both sides of \eqref{limit3} are finite, and then
\eqref{limit1}--\eqref{limit2} hold in the general case by dominated
convergence.

Thus,
$$\left(\int_{t\in\Om} [X'\cap
Y_t]\,\al(t)\,,\,\phi\right) =  \int_\Om \al(t)
\int _{X'\cap Y_t} \phi = \int_{\ycal'} \pi'^*_2\al \wedge
\pi'^*_1\phi = \big( [\ycal'] \wedge \pi_2^*\al,\,\pi_1^*\phi\big),$$
and therefore
\begin{equation}\label{limit4} \int_{t\in\Om} [X'\cap
Y_t]\,\al(t)= \pi_{1*} \big( [\ycal'] \wedge
\pi_2^*\al\big)\,.\end{equation}
Furthermore, since $\pi_1$ is a submersion, $\pi_1^*$ is well-defined on
currents and
\begin{equation}\label{limit5} \pi_{1*} \big( [\ycal'] \wedge
\pi_2^* \al\big) = \pi_{1*} \big( \pi_1^*[X'] \wedge
\pi_2^* \al\big) = [X']\wedge \pi_{1*}\pi_2^* \al\,.\end{equation}
The identity of the lemma follows from
\eqref{limit4}--\eqref{limit5} and \eqref{pushforward}. \end{proof}

We now proceed to the proof of Proposition \ref{converge}.  By abuse of notation, we
let $Z_{s^N}$ denote the zero set of a section $s^N$ as well as the current of
integration over the zero set. We choose
$N_0$ such that if
$N\ge N_0$, the zero sets $Z_{s^N_1},\dots, Z_{s^N_m}$ are almost always smooth
and intersect transversely.  (This holds if the Kodaira map for
$H^0(M,L^N)$ is an embedding.)

We begin by smoothing currents (locally) on $M$.
Let $a\in M$ be arbitrary and consider a coordinate chart
$\tau_a:V_a\buildrel \approx\over\to
B_r: = \{z\in\C^{m}:\|z\|<r\}$, with $a\in V_a\subset M,\ \tau(a)=0$. We let
$U_a= \tau_a\inv(B_{r/2})$.
To simplify our argument below, we  choose the biholomorphism
$\tau$ as follows: Embed
$M\subset \CP^q$, and choose
projective coordinates $(\zeta_0:\dots:\zeta_q)$ in $\CP^q$ such
that
\begin{itemize}\item $a=(1:0:\dots:0)$, \item $\{x\in M :
\zeta_j(x)=0\ \mbox{for }\ 0\le j\le m\} =\emptyset$,
\item the projection $\pi_a:M\to \CP^{m}$, $x\mapsto(\zeta_0(x)
:\dots:\zeta_{m}(x))$ has nonsingular Jacobian at $a$.
\end{itemize}  We  choose a neighborhood $V'_a$ of $a$
such that  $\pi_a$ is injective
on
$V'_a$ and $\pi_a(V_a')\subset \C^m=\CP^m\sm\{\zeta_0=0\}$. We
then choose
$r>0$ such that $B_r\subset \pi_a(V'_a)$, and we let
$V_a=\tau_a\inv(B_r)$ and  $\tau_a=\pi_a|_{V_a}$.

The advantage of this
construction is that degree bounds in $M$ push
forward under  $\tau_a$ to degree bounds in $\CP^{m}$.  In
particular, if $X$ is an algebraic hypersurface in $M$, then
\begin{equation}\label{push}\tau_a(X\cap V_a)\subset \pi_a(X),
\quad
\deg_{\CP^{m}}
\pi_{a}(X)
= \deg_{\CP^q}X\;.\end{equation}
(The well-known  formula \eqref{push} is easily
verified by recalling that the degree of a subvariety $X$ in projective
space is the number of points in the intersection of $X$ with a generic
linear subspace of complementary dimension.)

Let
$\psi_\ep(z)= \ep^{-2m}\psi(z/\ep)$ be an approximate identity on
$\C^m$, with $\psi\in\ccal^\infty(\C^m)$ and $\supp\psi\subset B_{r/2}$.  We
consider the local smoothing operator $S^a_\ep:\dcal'^{1,1}(M)\to \ecal^{1,1}(U_a)$
given by convolution in the
$\tau$ coordinates:
\begin{equation}\label{smoothZ} S^a_\ep u=   u *\psi_\ep:=
\tau_a^*\big[\tau_{a*}(u|_{V_a}) * \psi_\ep\big]\in
\ecal^{1,1}(U_a)\;, \quad\mbox{for }\ u\in \dcal'^{1,1}(M)\,.\end{equation}
(Note that $\tau_{a*}(u|_{V_a})\in\dcal'^{1,1}(B_r)$, and hence its
convolution with $\psi_\ep$ is well-defined on $B_{r/2}$ for $0<\ep< 1$.)

Now suppose that $P_0=(a,b)\in M\times M$, and let $\tau_a:V_a\buildrel
\approx\over\to B_r, \ \tau_b:V_b\buildrel \approx\over\to B_r$ be
as above, and let $\Om=U_a\times U_b$. We consider the approximate
identity $\wt\psi_\ep(z,w) = \psi_\ep(z)\psi_\ep(w)$ on $\C^{2m}$
and we similarly define $S_\ep:\dcal'^{2,2}(M)\to \ecal^{2,2}(\Om)$ by
\begin{equation}\label{smoothK} (S_\ep \wt u) = \wt u*\wt\psi_\ep:=
\tau^*\big[\tau_{*}(\wt u|_{V_a\times V_b}) * \wt\psi_\ep\big]\in
\ecal^{2,2}(\Om)\,, \quad\mbox{for }\ \wt u\in
\dcal'^{2,2}(M\times M)\,,\end{equation} where
$\tau=\tau_a\times\tau_b:V_a\times V_b\buildrel
\approx\over\to B_r \times B_r$.

\begin{lem} \label{smoothK21}For $N\ge N_0$, $$S_\ep K_{21}^N =
\E\big(S_\ep(Z_{s^N}
\boxtimes  Z_{s^N})\big)=
\E\big(S^a_\ep Z_{s^N}
\boxtimes S^b_\ep Z_{s^N} \big) \in \ecal^{2,2}(\Om)\,.$$ \end{lem}
\begin{proof} We have
$$ S_\ep K_{21}^N = \big[\E(Z_{s^N}
\boxtimes  Z_{s^N})\big]* \wt\psi_\ep \ = \
\E \big[(Z_{s^N}
\boxtimes  Z_{s^N})* \wt\psi_\ep\big] = \E\big(S_\ep(Z_{s^N}
\boxtimes  Z_{s^N})\big).$$
Furthermore,
$$\E \big[(Z_{s^N}
\boxtimes  Z_{s^N})* \wt\psi_\ep\big] =
 \E\big[(Z_{s^N}*\psi_\ep) \boxtimes (Z_{s^N}*\psi_\ep)\big]
\ =\ \E\big(S^a_\ep Z_{s^N}
\boxtimes S^b_\ep Z_{s^N} \big)\,.$$\end{proof}

\begin{lem}\label{intersect} Let $2\le k\le m$, $N\ge N_0$.
For almost all
$(s_1^N,\dots,s^N_k)\in H^0(M,L^N)^k$, we have
$$ [Z_{s^N_1} \boxtimes Z_{s^N_1}]\wedge
S_\ep[Z_{s^N_2} \boxtimes Z_{s^N_2}] \wedge \cdots\wedge
S_\ep[Z_{s^N_k} \boxtimes Z_{s^N_k}] \to
Z_{s^N_1,\dots,s^N_k}\boxtimes Z_{s^N_1,\dots,s^N_k}\,,$$
weakly in $\dcal'^{2k,2k}(\Om)$, as $\ep\to 0$.\end{lem}
\begin{proof}  It suffices to consider the case where the $Z_{s^N_j}$ are
smooth and intersect transversely. We let $Y_j =\big( Z_{s^N_j}
\times Z_{s^N_j}\big)\cap( V_a\times V_b)$, and we identify
$V_a\times V_b$ with $B_r\times B_r$ via the biholomorphism
$\tau$. Under this identification, $\Om=B_{r/2}\times
B_{r/2}\subset \C^{2m}$. For $t\in \C^{2m}$, let
$T_t:\C^{2m}\to\C^{2m}$ denote the translation $T_t(w)=w+t$,
so that $$S_\ep Y_j = \int_{t\in\Om}
\left[T_{-t}\,Y_j\right] \wt\psi_\ep(t)\,\nu(t)\in
\ecal^{2,2}(\Om)\,,
$$ where $\nu$ is the Euclidean volume form on $\C^{2m}$.

Suppose that $X'$ is a complex submanifold of $\Om$. We first show by
induction that
$X'
\cap T_{t_2}Y_2
\cap
\cdots\cap T_{t_k}Y_k$ is a complex submanifold of $\Om$ for almost
all
$t_2,\dots, t_k$, and
\begin{multline}\label{Xwedge}
[X']\wedge S_\ep Y_2\wedge\cdots\wedge S_\ep Y_k
\\= \int_{\Om^{k-1}}\left[X' \cap T_{t_2}Y_2 \cap \cdots\cap
T_{t_k}Y_k\right] \wt\psi_\ep(t_2)\cdots \wt
\psi_\ep(t_k)\,\nu(t_2)\wedge\cdots\wedge \nu(t_k).\end{multline}
To verify \eqref{Xwedge} for $k=2$, we let $\ycal=\{(z,t)\in\Om\times \Om:
z-t\in Y_2\}$.  Then $\ycal$ is smooth and the two projections
$\pi_1:\ycal\to \Om$, $\pi_2:\ycal\to \Om$ are
submersions.  Furthermore, $\pi_1(\pi_2\inv\{t\})= T_{t}Y_2$. Hence
by Lemma \ref{intersectEZ} with $X=\Om$ and
$\al=\wt\psi_\ep\nu$, the intersection $X' \cap T_{t_2}Y_2$ is
a complex submanifold for almost all $t_2\in\Om$, and
\eqref{Xwedge} holds for $k=2$. For the inductive step, let $k>2$ and
suppose that
 \eqref{Xwedge} has been verified
for $k-1$. Let $t_1,\dots,t_{k-1}$ be parameters in $\Om$ such
that  $X'\cap T_{t_2}Y_2\cap
\cdots\cap T_{t_{k-1}}Y_{k-1}$ is a complex submanifold of $\Om$. By
Lemma \ref{intersectEZ} with $X'$ replaced by $X'\cap T_{t_2}Y_2\cap
\cdots\cap T_{t_{k-1}}Y_{k-1}$ and $\ycal=\{(z,t)\in\Om\times \Om:
z-t\in Y_k\}$, we conclude that $X'\cap T_{t_2}Y_2\cap
\cdots\cap T_{t_k}Y_k$ is a complex submanifold of $\Om$ for almost
all $t_k$, and
$$\int_\Om \left[X'\cap Y^1_{t_1}\cap \dots\cap
Y^{k-1}_{t_{k-1}}\cap Y^k_{t_k}\right] \psi_\ep(t_k)\,\nu(t_k) =
\left[X'\cap Y^1_{t_1}\cap
\dots\cap Y^{k-1}_{t_{k-1}}\right] \wedge S_\ep Y_k\;.$$ Integrating
over
$t_1,\dots,t_{k-1}$ and applying the inductive assumption, we obtain
\eqref{Xwedge}.

Setting $X'=Y_1$ in \eqref{Xwedge}, we have
\begin{multline}\label{Ywedge}
\big[Y_1\big]\wedge S_\ep Y_2\wedge\cdots\wedge S_\ep Y_k
\\= \int_{\Om^{k-1}}\left[Y_1 \cap T_{t_2}Y_2 \cap \cdots\cap
T_{t_k}Y_k\right] \wt\psi_\ep(t_2)\cdots \wt
\psi_\ep(t_k)\,\nu(t_2)\wedge\cdots\wedge \nu(t_k).\end{multline}
Now choose $\ep_0>0$ such that $Y_1 , T_{t_2}Y_2 ,\dots,
T_{t_k}Y_k$ intersect transversely whenever $|t_j|<\ep_0$ for $2\le
j\le k$.
Let $\phi\in\dcal^{2m-2k,2m-2k}(\Om)$ be a test form.  Since the
submanifolds $Y_1 \cap T_{t_2}Y_2 \cap\cdots\cap
T_{t_k}Y_k$ vary smoothly as the parameters $t_2,\dots,t_k$ vary in
the $\ep_0$-ball, it follows from an argument using the implicit function
theorem that the map
$$(t_2,\dots,t_k)\mapsto
\int_{Y_1 \cap T_{t_2}Y_2 \cap\cdots\cap
T_{t_k}Y_k} \phi=
\left(\left[Y_1 \cap T_{t_2}Y_2 \cap \cdots\cap
T_{t_k}Y_k\right],\,\phi\right)$$ is continuous (and in fact is
$\ccal^\infty$) for
$|t_j|<\ep_0$.  Therefore by
\eqref{Ywedge},
\begin{equation*}  \left(\big[Y_1\big]\wedge S_\ep
Y_2\wedge\cdots\wedge S_\ep Y_k,\phi\right)\to
\left(\left[Y_1 \cap Y_2 \cap \cdots\cap
Y_k\right],\,\phi\right) \qquad \mbox {as }\
\ep\to 0\;;\end{equation*} i.e.,
\begin{equation}
\label{weakk}\big[Y_1\big]\wedge S_\ep
Y_2\wedge\cdots\wedge S_\ep Y_k\to \left[Y_1
\cap Y_2 \cap \cdots\cap Y_k\right] \qquad \mbox {weakly, \ \
as }\ \ep\to 0\;.\end{equation}\end{proof}

\begin{lem} \label{unifbound}There exists a
positive constant $C<+\infty$ such that for all $N\ge N_0$,
$\phi\in\dcal^{2m-2k,2m-2k} (\Om)$, and $0<\ep<1$, we
have
$$\left|[Z_{s^N_1} \boxtimes Z_{s^N_1}]\wedge
S_\ep[Z_{s^N_2} \boxtimes Z_{s^N_2}] \wedge \cdots\wedge
S_\ep[Z_{s^N_k} \boxtimes Z_{s^N_k}],\phi\big)\right| \le C
N^{2k} \|\phi\|_\infty$$ for almost all $(s_1^N,\dots s_k^N)\in
H^0(M,L^N)^k$.
\end{lem}

\begin{proof} Fix $s_1,\dots, s_k$ so that $Z_{s^N_1},\dots,Z_{s^N_k}$
intersect transversely, and let $Y_j=\big( Z_{s^N_j}
\times Z_{s^N_j}\big)\cap( V_a\times V_b)$, as in the proof of Lemma \ref{intersect}.
By
\eqref{Ywedge} it suffices to show that
\begin{equation}\label{bound1} \left(\left[Y_1 \cap T_{t_2}Y_2
\cap \cdots\cap
T_{t_k}Y_k\right],\,\phi\right) \le C
N^{2k} \|\phi\|_\infty\;,\end{equation} for almost all $(t_2,\dots,t_k)
\in \Om^{k-1}$. To verify \eqref{bound1}, it in turn suffices to show that
\begin{equation}\label{bound2} \vol_{M\times M}\left(Y_1 \cap
T_{t_2}Y_2
\cap \cdots\cap
T_{t_k}Y_k\right) \le C
N^{2k}\;.\end{equation} We write
$t_j=(t_j',t_j'')$, $Z^a_j=Z_{s^N_j}\cap V_a$, $Z^b_j=Z_{s^N_j}\cap V_b$,
so that $T_{t_j} Y_j=T_{t'_j}Z^a_j \times T_{t''_j}Z^b_j$.  Then
\begin{multline}\label{bound3} \vol_{M\times M}\left(Y_1 \cap
T_{t_2}Y_2
\cap \cdots\cap
T_{t_k}Y_k\right) \\= \vol_M (Z^a_1 \cap  T_{t'_2} Z^a_2
\cap\cdots\cap T_{t'_k}Z^a_k)\;  \vol_M
(Z^b_1 \cap  T_{t''_2} Z^b_2
\cap\cdots\cap T_{t''_k}Z^b_k)\;.\end{multline}
Since $\tau_a:Z^a_j
\hookrightarrow
\pi_a(Z_{s^N_j})\subset \CP^{m}$ is injective and  the translations
$T_{t_j'}$ extend to automorphisms of $\CP^{m}$, we have
\begin{eqnarray*}\vol_M\left(Z^a_1 \cap  T_{t'_2} Z^a_2
\cap\cdots\cap T_{t'_k}Z^a_k\right) &\le& C_1\,
\vol_{\CP^m}\left(\pi_a(Z_{s^N_1})\cap T_{t'_2}\pi_a(Z_{s^N_2})
\cap \cdots\cap
T_{t'_k}\pi_a(Z_{s^N_k})\right)\\&=&\frac{\pi^kC_1}{k!}
\,\prod_{j=1}^k \deg_{\CP^m} \pi_a(Z_{s^N_j})\;.
\end{eqnarray*}
However, by \eqref{push}, \begin{equation*} \deg_{\CP^m}
\pi_a(Z_{s^N_j}) = \deg_{\CP^q} Z_{s^N_j} = \frac {1}{\pi^{m-1}}
\int_{Z_{s^N_j}}  \om_{\CP{^q}}^{m-1} = \frac {1}{\pi^{m-1}}
\int_M N\,c_1(L,h)\wedge \om_{\CP{^q}}^{m-1}=C_2N,\end{equation*}
and hence
\begin{equation}\label{bound4}\vol_M\left(Z^a_1 \cap  T_{t'_2} Z^a_2
\cap\cdots\cap T_{t'_k}Z^a_k\right) \le C_3 N^k\;.\end{equation}
The bound \eqref{bound2} follows from
\eqref{bound3}--\eqref{bound4}.
\end{proof}

\begin{lem} \label{easy} Let  $f_j\in
L^n_{loc}(\R^k)$ for $1\le j\le n$. Let $f_j^\ep= f_j*\psi_\ep$, where
$\psi_\ep$ is a compactly supported smooth approximate identity. Then
$$\prod_{j=0}^n f_j^\ep\ \buildrel{L^1_{loc}}\over\longrightarrow\
\prod_{j=0}^n f_j\qquad
\mbox{as }\ \ep\to 0\;.$$
\end{lem}

\begin{proof} We use the
generalized H\"older inequality:
\begin{equation}\label{Holder}\sum_{j=1}^n\frac 1{p_j}=1\
\implies\
\|f_1\cdots f_n\|_1 \le \|f_1\|_{p_1}
\cdots
\|f_n\|_{p_n}\;.\end{equation}
We can assume without loss of
generality that the
$f_j$ have compact support and hence  $f_j\in
L^n(\R^k)$, for $1\le j\le n$.  By
\eqref{Holder} with
$p_j=n$, we then have
\begin{eqnarray*}\| f_1^\ep \cdots f_n^\ep - f_1\cdots f_n\|_1
\hspace{-1.5in}\\&\le&
 \|(f_1^\ep - f_1) f_2^\ep
\cdots f_n^\ep\|_1\ + \|f_1(f_2^\ep-f_2)f_3^\ep\cdots f_n^\ep\|_1+ \cdots+\
\|f_1\cdots f_{n-1} (f_n^\ep-f_n)\|_1\\&\le&
\|(f_1^\ep - f_1)\|_n\,\| f_2^\ep\|_n\,
\cdots \|f_n^\ep\|_n\  +\ \cdots\ +\
\|f_1\|_n\,\cdots \|f_{n-1}\|_n\, \|(f_n^\ep-f_n)\|_n\
\to\  0\;.\end{eqnarray*}
\end{proof}

Part (i) of Proposition \ref{converge} is an immediate consequence of
Lemma \ref{smoothK21}. Next we show part (iii):  By Lemma \ref{smoothK21}
and the independence of the
$s^N_j$, we have \begin{eqnarray*} \K^N_{21} \wedge
\left(S_\ep \K^N_{21}\right)^{k-1}  &=&
 \E\big(Z_{s^N_1} \boxtimes
Z_{s^N_1}\big)\wedge
\E\big(S_\ep (Z_{s^N_2} \boxtimes  Z_{s^N_2})\big) \wedge \cdots\wedge
\E\big(S_\ep (Z_{s^N_k} \boxtimes  Z_{s^N_k})\big)\\
&=& \E\big([Z_{s^N_1} \boxtimes
Z_{s^N_1}]\wedge S_\ep (Z_{s^N_2} \boxtimes  Z_{s^N_2}) \wedge
\cdots\wedge S_\ep (Z_{s^N_k} \boxtimes  Z_{s^N_k})\big)\,.\end{eqnarray*}

Therefore, for a test form $\phi\in\dcal^{m-k,m-k}(\Om)$, we
have
\begin{multline} \left( \K^N_{21}
\wedge \left(S_\ep \K^N_{21}\right)^{k-1}, \phi\right)\\=
\int_{H^0(M,L^N)^k} \left([Z_{s^N_1} \boxtimes
Z_{s^N_1}]\wedge S_\ep (Z_{s^N_2} \boxtimes  Z_{s^N_2}) \wedge
\cdots\wedge S_\ep (Z_{s^N_k} \boxtimes  Z_{s^N_k}),\phi\right)\left[
\prod_{j=1}^k d\ga_N(s_j^N)\right].\label{sm4}\end{multline} By
Lemma \ref{unifbound}, the integrand in \eqref{sm4} is uniformly
bounded, and hence by \eqref{pcc}, Lemma \ref{intersect} and
Lebesgue dominated convergence, we have
$$ \left( \K^N_{21}
\wedge \left(S_\ep \K^N_{21}\right)^{k-1}, \phi\right)\to
\int_{H^0(M,L^N)^k} \left(Z_{s^N_1,\dots,s^N_k}\boxtimes
Z_{s^N_1,\dots,s^N_k},\phi\right)\left[ \prod_{j=1}^k
d\ga_N(s_j^N)\right]= \left(\K_{2k}^N,\phi\right),$$ as $\ep\to
0$, verifying part (iii).

To complete the proof of the proposition, we recall from
Proposition
\ref{KN21} that the current
$\K_{21}^N$ has $L^{m-1}_{loc}$ coefficients (if $m\ge 2$) and hence by
Lemma \ref{easy},
\begin{equation}\label{iii-} (S_\ep \K_{21}^N)^j\to (\K_{21}^N)^j|_\Om
\quad \mbox {in }\ L^1(\Om), \quad\mbox{for }\ 1\le j\le
m-1\,.\end{equation} Since $L^1$ convergence implies weak
convergence, it follows from (iii) and \eqref{iii-} that
$\K^N_{2j}|_\Om = (\K_{21}^N)^j|_\Om$ and hence (ii) holds. This
completes the proof of Proposition \ref{converge}.\qed

\subsection{Completion of the proof of Theorem \ref{varint2}}
We recall that
\begin{equation}\label{recall}\vol_{2m-2k}[Z_{s_1^N,\dots,s_k^N}\cap U] =
(Z_{s_1^N,\dots,s_k^N},\chi_U\Phi_k)=(Z_{s_1^N,\dots,s_k^N},\chi_{\overline
U} \,\Phi_k)\qquad a.s.\end{equation} (To verify the second equality in
\eqref{recall}, we note that $(Z_{s_1^N,\dots,s_k^N}, \chi_{\d
U}\Phi_k)= 0$
almost surely, since $\E\big(Z_{s_1^N,\dots,s_k^N}, \chi_{\d
U}\Phi_k\big)=0$ by Corollary \ref{indep}.)
We now approximate $\chi_{\overline U}$ by a sequence of
$\ccal^\infty$ functions $\chi_n:M\to \R$, $n=1,2,3,\dots,$
satisfying:
\begin{itemize}\item $0\le \chi_n\le 1$,
\item $\sup|d\chi_n| =O(n)$,\item $\chi_n|_{\overline U}\equiv 1$,\item
$\chi_n(w)=0$ for dist$(U,w)> 1/n$. \end{itemize} To construct
$\chi_n$, we choose $\rho\in\ccal^\infty(\R)$ such that
$\rho(t)=1$ for $t\le \frac 13$, $\rho(t)=0$ for $t\ge \frac 23$,
and $0\le \rho\le 1$. Let $\chi^0_n(w)=\rho(n\,\mbox{dist}(U,w))$.
If $\d U$ is smooth, then $\chi_n^0$ is smooth, for $n$
sufficiently large, and we can take $\chi_n=\chi_n^0$. Otherwise,
the Lipschitz constant of $\chi_n^0$ is $O(n)$, and we can smooth
$\chi_n^0$ to obtain our desired $\ccal^\infty$ function $\chi_n$.

Then $\chi_n \to \chi_{\overline U}$ pointwise, and hence for all
$(s^N_1,\dots,s^N_k)$, we have by Lebesgue dominated convergence,
$$(Z_{s_1^N,\dots,s_k^N},\chi_n\Phi_k) \to
(Z_{s_1^N,\dots,s_k^N},\chi_{\overline U}\,\Phi_k) =
\vol_{2m-2k}[Z_{s_1^N,\dots,s_k^N}\cap \overline U] \qquad
\mbox{as }\ n\to\infty\,.$$ Therefore (again by  dominated
convergence),
\begin{equation}\label{lim}
\var\big(Z_{s_1^N,\dots,s_k^N},\chi_n\Phi_k\big)\to \var
\big(\vol_{2m-2k}[Z_{s_1^N,\dots,s_k^N}\cap \overline U]\big) =
\var \big(\vol_{2m-2k}[Z_{s_1^N,\dots,s_k^N}\cap U]\big)
\end{equation} as $ n\to\infty$. To complete the proof of
Theorem \ref{varint2}, it suffices by \eqref{lim} and  Theorem \ref{variant}
with
$\phi=\chi_n\Phi_k$ to show that
\begin{multline} \label{lim2}\int_{M\times M}\dbar_1\dbar_2
Q_N\wedge\left(\d_1\dbar_1 \d_2\dbar_2 Q_{N}\right)^{j-1}\wedge
\big(\E Z_{s^N} \boxtimes \E Z_{s^N}\big)^{k-j}\wedge \big(
\d[\chi_n\Phi_k]\boxtimes \d[\chi_n\Phi_k]\big)\\= \int_{M\times
M}\dbar_1\dbar_2 Q_N\wedge\left(\d_1\dbar_1 \d_2\dbar_2
Q_{N}\right)^{j-1}\wedge \big(\E Z_{s^N} \boxtimes \E
Z_{s^N}\big)^{k-j}\wedge(\Phi_k\boxtimes\Phi_k)\wedge
(d\chi_n\boxtimes d\chi_n)\\
\to  - \int_{\d U\times \d U}  \dbar_1\dbar_2 Q_N\wedge\left(\d_1\dbar_1
\d_2\dbar_2 Q_{N}\right)^{j-1}\wedge \big(\E Z_{s^N} \boxtimes \E
Z_{s^N}\big)^{k-j}\wedge(\Phi_k\boxtimes\Phi_k).\quad\end{multline}

To verify \eqref{lim2}, let $$f= \dbar_1\dbar_2Q_N\wedge
\left(\d_1\dbar_1 \d_2\dbar_2 Q_{N}\right)^{j-1}\wedge \big(\E
Z_{s^N} \boxtimes \E Z_{s^N}\big)^{k-j}\wedge
(\Phi_k\boxtimes\Phi_k)\;.$$ We must show that $f|_{\d U\times \d
U}$ is $L^1$ and
\begin{equation}\label{lim3} \int_{M\times M} f\wedge (d\chi_n\boxtimes
d\chi_n)\to   - \int_{\d U\times \d U} f\;.  \end{equation}

 By Lemma
\ref{DISTM2}, we have \begin{equation}\label{f}|f(z,w)| =
O\left(N^j\,\mbox{dist}(z,w)^{-2j+2}\right)\le
O\left(N^j\,\mbox{dist}(z,w)^{-2m+2}\right)\;.\end{equation} Since
$\d U$ is a finite union of $\ccal^2$ submanifolds of $M$ of real
dimension $2m-1$, it follows from \eqref{f} that $f$ is $L^1$ on
$\d U\times \d U$.

Let $\de>0$ and consider the cut-off function $\la_\de(z,w)=
\rho(\de\inv \mbox{dist}(z,w))$, where $\rho$ is as above.  Then
$\la_\de\in\ccal^\infty(M\times M)$ for $\de$ sufficiently small,
$\la_\de(z,w)=0$ if dist$(z,w)>\de$, and $\la_\de(z,w)=1$ if
dist$(z,w)<\de/3$. We decompose the integral in \eqref{lim3}:
\begin{equation}\label{decomp} \int_{M\times M} f\wedge (d\chi_n\boxtimes
d\chi_n) = \int_{M\times M}\la_\de\, f\wedge (d\chi_n\boxtimes
d\chi_n)+\int_{M\times M} (1-\la_\de)\,f\wedge (d\chi_n\boxtimes
d\chi_n)\,.\end{equation} Since $(1-\la_\de)\,f$ is smooth and
$\chi_n\to\chi_U$,   it follows  that
\begin{multline}\label{lim4}\int_{M\times M} (1-\la_\de)\,f\wedge
(d\chi_n\boxtimes d\chi_n) = (d\chi_n\boxtimes d\chi_n,
(1-\la_\de)\,f)\\ \to (d\chi_U\boxtimes d\chi_U,
(1-\la_\de)\,f)=-\int_{\d U\times \d U}
(1-\la_\de)\,f\;.\end{multline} 
(The minus sign in \eqref{lim3} is due to the fact that
$(A\boxtimes B,
\phi\boxtimes \psi) = (-1)^{\deg B\,\deg \phi}(A,\phi)(B,\psi)$, and hence 
$d\chi_U\boxtimes d\chi_U = [-\d U]\boxtimes [-\d U]= -[\d U\times \d U]$,
where $[\d U]$ denotes the current of integration over $\d U$.) 

To complete the proof   
of \eqref{lim3},  we must show that the
$\la_\de\,f$ integrals are uniformly small.  For $z_0\in M,\ n\in\Z^+,\ \de>0$,
we write
\begin{eqnarray*}V(z_0,n,\de) &:=&\{w\in M: \mbox{dist}(z,w)< \de,\ w\in
\supp(d\chi_n)\}\\& \subset&\{w:\mbox{dist}(z,w)< \de,\
\mbox{dist}(U,w)<1/n \}\,.\end{eqnarray*}  Since $\d U$ is
piecewise smooth, we can choose $\de_0>0,\ n_0\in\Z^+$ such that
for all $z_0\in M$: \begin{itemize} \item the exponential map
$\exp_{z_0}:T_{z_0}(M) \to M$ is injective on the $\delta_0$-ball
$B_{\de_0}(z_0)=\{v\in T_{z_0}(M):|v|<\de_0\}$; \item there exists
real hyperplanes $P_1,\dots,P_q$, such that
\begin{equation}\label{platter}V(z_0,n,\de_0)
\subset\bigcup_{j=1}^q\exp_{z_0}\big(\{v+tu_j\in B_{\de_0}(z_0):
v\in P_j,\ |t_j|<2/n\}\big)\,,\end{equation} for all $n>n_0$,
where $u_j$ is a unit normal to $P_j$.\end{itemize} Here, $q$ is
the maximal number of facets of the polyhedral cones locally
diffeomorphic to open sets of $\d U$, as described after the
statement of Theorem \ref{number}. (If $\d U$ is smooth, then
$q=1$.)

Since $j\le m$ and $|d\chi_n|=O(1/n)$, we then have by \eqref{f}
and \eqref{platter},
\begin{eqnarray*}\left| \int_M \la_\de(z_0,w)f(z_0,w)\wedge
d\chi_n(w)\right| &\le & \left|\int_{\{z_0\}\times V(z_0,n,\de)}
f(z_0,w)\wedge
d\chi_n(w) \right| \\
&\le & Cn \int_{V(z_0,n,\de)} \mbox{dist}(z_0,w)^{-2m+2}\,d\vol_M\\
&\le & C'n\int_{\{x\in \R^{2m}:|x|<\de,\ |x_1|<2/n\}}
|x|^{-2m+2}\,dx\\ &\le & 4C'\int_{\{y\in \R^{2m-1}:|y|<\de\}}
|y|^{-2m+2}\,dy\ =\ C''\de\,,
\end{eqnarray*} where $C,C',C''$ are constants independent of $z_0$
(but depending on $m,U,N$).  Here, $f(z_0,w)$ is regarded as a
$(2m-1)$-form (in the $w$ variable) with values in
$T^{*2m-1}_{z_0}(M)$.  Therefore,
\begin{eqnarray*}\left|\int_{M\times M}\la_\de\, f\wedge
(d\chi_n\boxtimes d\chi_n)\right|&=&\left|\int_{\{z\in M:
\operatorname{dist}(U,z)<1/m\}}d\chi_n(z)\int_{\{z\}\times
M}\la_\de(z,w)\, f(z,w)\wedge d\chi_n(w)\right|\\&\le & C''\de
\int_{\{z\in M: \operatorname{dist}(U,z)<1/m\}}
|d\chi_n(z)|\,d\vol_{\d U}(z)\\&\le & C''\de\, \sup |d\chi_n|\,
\vol(\{z\in M: \operatorname{dist}(U,z)<1/n\})\,.
\end{eqnarray*}  Since $\sup |d\chi_n|=O(n)$ and the volume of the
shell $\{z\in M: \operatorname{dist}(U,z)<1/n\}$ is $O(1/n)$, it
follows that \begin{equation}\label{db}\left|\int_{M\times
M}\la_\de\, f\wedge (d\chi_n\boxtimes d\chi_n)\right|\le C'''\de
\qquad \forall\ n>n_0\,.
\end{equation} Then \eqref{lim3} follows from \eqref{decomp},
\eqref{lim4} and \eqref{db}, which completes the proof of Theorem
\ref{varint2}.\qed

\medskip
\section{Variance of zeros in a domain:
Proof of Theorems \ref{number} and \ref{volume}}\label{s-number}

We now use Theorem \ref{varint2} together with the asymptotics of the
pluri-bipotential $Q_N$ to prove Theorem \ref{volume}.

By Theorem \ref{varint2} and  Proposition \ref{indepTYZ}, we have
\begin{equation}\label{VU}\var\big(\vol_{2m-2k}[Z_{s_1^N,\dots,s_k^N}\cap U]\big)
=\sum_{j=1}^k {k\choose j}\,V_j^N(U)\;,\end{equation} where
\begin{eqnarray}\label{VjU} V_j^N(U)&=&
\left(\frac N\pi \right)^{2k-2j}
\int_{\d U\times\d U}
  -\dbar_1\dbar_2  Q_N(z,w)\wedge\big[-
\d_1\dbar_1 \d_2\dbar_2
Q_{N}(z,w)\big]^{j-1}\nonumber\\
&&\qquad\wedge\left[\om(z)^{k-j}+O\left(\frac
1N\right)\right]\wedge\left[\om(w)^{k-j}+O\left(\frac
1N\right)\right]\wedge \Phi_k(z)\wedge\Phi_k(w)\nonumber\\&=&\frac
1{(m-k)!^2}\left(\frac N\pi \right)^{2k-2j} \int_{\d U}
\Upsilon^N_j\wedge\left[\om^{m-j} +O\left(\frac
1N\right)\right]\,,\end{eqnarray} where $\Phi_k$ is given by
\eqref{Phi}, and
\begin{eqnarray}\label{int2}\Upsilon^N_j(z)&:=& \int_{\{z\}\times \d
U}  - \dbar_1\dbar_2  Q_N(z,w)\wedge\big[- \d_1\dbar_1 \d_2\dbar_2
Q_{N}(z,w)\big]^{j-1}\wedge\left[\om(w)^{m-j}+O\left(\frac
1N\right)\right]\nonumber\\&&\qquad\in
T^{*j-1,j}_z(M)\;.\end{eqnarray}
  (In \eqref{int2} and below, we regard the integrand as an $(m-1,m)$-form
in the $w$ variable with values in $T^{*j-1,j}_z(M)$ by identifying 
$\pi_1^*\phi \wedge \pi_2^*\psi\in T_{z,w}^{*m+j-2,m+j}(M\times
M)$ with 
$\phi\otimes \psi\in T^{*j-1,j}_z(M)\otimes T^{*m-1,m}_w(M)$, for $\phi\in
T^{*j-1,j}_z(M)$, $\psi\in T^{*m-1,m}_w(M)\,$.)
 
By Lemma \ref{Qdecay},
\begin{equation}\label{d2Qdecay}\textstyle
\dbar_1\dbar_2Q_N(z,w)\wedge\big[ \d_1\dbar_1 \d_2\dbar_2
Q_{N}(z,w)\big]^{j-1}= O(N^{-m})\;,
 \quad \mbox{for }\ \dist(z,w)>b\sqrt{\frac{\log N}N}\;,\end{equation}
where we choose $b=\sqrt{2m+3}$.  Thus we can approximate
$\Upsilon^N_j(z)$ by restricting the integration in \eqref{int2}
to the set of $w\in\d U$ with $\dist(z,w)<b\sqrt {\frac{\log N}N}$.

To evaluate $\Upsilon^N_j(z_0)$ at a fixed point $z_0\in\d U$, we
choose normal holomorphic coordinates $\{w_1,\dots,w_m\}$ centered
at $z_0$ and defined in a neighborhood $V$ of $z_0$, and we make
the change of variables $w_j=\frac {v_j}\sqrtn$  as in \S \ref{off}.
 Since $\om=\frac i2 \ddbar \log a= \frac i2 \ddbar \left[|w|^2
+O(|w|^3)\right]$, we note that
\begin{equation}\label{kahlerE}
\om\left(z_0+\frac v{\sqrtn}\right)= \frac
i{2}\sum\left[\de_{jk}+O\left(\frac
{|v|}{\sqrtn}\right)\right]\frac 1N dv_j\wedge d \bar v_k= \frac
i{2N}\ddbar|v|^2+O\left(\frac{|v|}{N^{3/2}}\right)\;,\end{equation}
for $|v|\le b\sqrt{\log N}$.  We then have
\begin{eqnarray}\label{int3}\Upsilon^N_j(z_0)&=& N^{j-m}
\int_{\left\{|v|\le b\sqrt{\log N}:\,z_0+\frac v\sqrtn\in\d
U\right\}}\textstyle -\dbar_1\dbar_2 Q_N(z_0,z_0+\frac
v\sqrtn)\nonumber\\&& \textstyle\quad \wedge\big[- \d_1\dbar_1
\d_2\dbar_2 Q_{N}(z_0,z_0+\frac v\sqrtn)\big]^{j-1}
\wedge\left[(\frac i2\ddbar |v|^2)^{m-j}+O\left(\frac
1N\right)\right]\;.\end{eqnarray} Applying the asymptotics of
Lemmas \ref{d2Qas}--\ref{d4Qas} to \eqref{int3}, we obtain the
formula
\begin{eqnarray} \label{Upj}\Upsilon^N_j(z_0) &=& N^{2j-m-1/2}
\left[\int_{\left\{|v|\le b\sqrt{\log N}:\,z_0+\frac v\sqrtn\in\d
U\right\}}\textstyle\frac 14\,  F''(\half|v|^2)\,   \dbar(\bar z\cdot v)
\wedge \dbar|v|^2  \right.\nonumber\\&&
\left.\qquad\wedge
\left(\Var^{z_0}_\infty\right)^{j-1}\wedge\left(\frac i2\ddbar
|v|^2\right)^{m-j}\ +\ O(N^{-1/2+\ep})\right]\;.\end{eqnarray}

We first consider the case where $\d U$ is $\ccal^2$ smooth
(without corners). We can choose our holomorphic normal
coordinates $\{w_j\}$ so that the real hyperplane $\{\Im w_1=0\}$
is tangent to $\d U$ at $z_0$.   We can then write (after
shrinking the neighborhood $V$ if necessary),
$$U\cap V=\{w\in V:\Im w_1+\phi(w)>0\}\;,$$ where  $\phi:V\to\R$ is a
$\ccal^2$ function of $(\Re w_1,\, w_2,\dots,w_m)$ such that
$\phi(0)=0,\ d\phi(0)=0$.

We consider the {\it nonholomorphic\/} variables
\begin{equation}\label{wtilde} \wt w=\tau(w):=(w_1+i\phi(w),w_2,\dots,
w_m)\;,\end{equation} so that $\d U=\{\Im \wt w_1=0\}$. We next
make the change of variables
$$\wt v=\tau_N(v):=\sqrtn \,\tau\left(\frac v{\sqrtn}\right)=\sqrtn\,\wt w = v
\left[1+O\left(\frac v{\sqrtn}\right)\right]$$ in the integral
\eqref{Upj} to obtain
\begin{eqnarray}\label{univ0} \Upsilon^N_j(z_0) &= & N^{2j-m-1/2}\Big[
\int_{\{\wt v\in B_N^{2m-1}\}}\textstyle\frac 14\, F''(\half|\wt
v|^2)\,   \dbar(\bar z\cdot \wt v)
\wedge \dbar|\wt v|^2 \nonumber
\\&&\quad\textstyle \wedge \left(\Var^{z_0}_\infty(\wt
v)\right)^{j-1}\wedge\left(\frac i2\ddbar |\wt v|^2\right)^{m-j}\
+\ O(N^{-1/2+\ep})\Big],\end{eqnarray} where
$$\left\{v\in\R\times\C^{m-1}:|v|<(b-1)\sqrt{\log N}\right\}
\subset B_N^{2m-1}\subset
\left\{v\in\R\times\C^{m-1}:|v|<(b+1)\sqrt{\log N}\right\}.$$

By \eqref{F''} and
\eqref{varest}, we have $F''(\half |v|^2)|v|^2 \{
\Var^{z_0}_\infty(v)\}^{j-1} =O(e^{-|v|^2})$ for $|v|>1$, and
hence
\begin{multline*}\int_{|\wt v|>(b-1)\sqrt{\log N}}\textstyle\left|
F''(\half|\wt v|^2)\,   \dbar(\bar z\cdot\wt v)
\wedge \dbar|\wt v|^2  \wedge \left(\Var^{z_0}_\infty(\wt
v)\right)^{j-1}\wedge \left(\frac i2\ddbar |\wt
v|^2\right)^{m-j}\right|\\=O\left(N^{-(b-1)^2+\ep}\right) =
O\left(N\inv\right)\;.$$
\end{multline*}
Thus we can  replace the $B_N^{2m-1}$ integral in
\eqref{univ0} with the affine integral over $\R\times \C^{m-1}$,
so that
\begin{equation}\label{univ}
\Upsilon^N_j(z_0) = N^{2j-m-1/2}\left[\Upsilon^\infty_j(z_0) +
O(N^{-1/2+\ep})\right],\end{equation} where
\begin{eqnarray}\label{Upj1}\Upsilon^\infty_j(z_0) &\!\!:\,=& \int_{\R\times
\C^{m-1}}\textstyle\frac 14\, F''(\half|v|^2)\,   \dbar(\bar z\cdot v)
\wedge \dbar|v|^2  \wedge
\left(\Var^{z_0}_\infty\right)^{j-1}\wedge\left(\frac i2\ddbar
|v|^2\right)^{m-j}\nonumber\\&=& \frac 1{4\pi^2}\int_{\R\times
\C^{m-1}}\frac 1{e^{|v|^2}-1}\,   \dbar(\bar z\cdot v)
\wedge \dbar|v|^2 \wedge
\left(\Var^{z_0}_\infty\right)^{j-1}\wedge\left(\textstyle\frac
i2\ddbar |v|^2\right)^{m-j}.\end{eqnarray} 
  Since  only the last two terms of $\Var^{z_0}_\infty$ contain neither 
$\dbar(\bar z\cdot v)$ nor $ \dbar|v|^2$, formula \eqref{Upj1} simplifies to:
\begin{multline}\label{Upj1s}\Upsilon^\infty_j(z_0) = \frac
{(-1)^{j-1}}{(4\pi^2)^j}\int_{\R\times
\C^{m-1}}\frac 1{(e^{|v|^2}-1)^j}\, \dbar(\bar z\cdot v)
\wedge \dbar|v|^2 \\ \wedge
\big[ \dbar \d (\bar z\cdot v)\wedge\ddbar(z\cdot\bar v) +
\ddbar|z|^2\wedge \ddbar |v|^2\big]^{j-1}\wedge\left(\textstyle\frac
i2\ddbar |v|^2\right)^{m-j}.\end{multline} 
Thus,  \begin{equation} (\Upsilon^\infty_j\wedge \om^{m-j})(z_0) = c_{mj}
\;dx_1\wedge \textstyle (\frac i2 dz_2\wedge d\bar
z_2)\wedge\cdots\wedge (\frac i2 dz_m\wedge d\bar z_m)=
c_{mj}\;d\vol_{\d U,z_0}\,,\label{Upj2}\end{equation} where
$c_{mj}$ is a universal constant,   which we compute in \S \ref{positivity}
below.  

Substituting \eqref{univ} and  \eqref{Upj2} in \eqref{VjU}, we
have \begin{equation}\label{VjU1} V_j^N(U)= \frac
1{(m-k)!^2}\left(\frac 1{\pi^{2k-2j}} \right)N^{2k-m-1/2}\left[
\int_{\d U} c_{mj}\;d\vol_{\d U,z_0}
+O\left(N^{-1/2+\ep}\right)\right].\end{equation} Combining
\eqref{VU} and \eqref{VjU1}, we obtain the formula of Theorem
\ref{volume} with \begin{equation}\label{numk}\nu_{mk} = \frac
1{(m-k)!^2}\sum_{j=1}^k {k\choose j} \frac
{c_{mj}}{\pi^{2k-2j}}\,,\end{equation} for the case where $\d U$
is smooth.

  We  now  verify the general case where $\d U$ is piecewise
smooth (without cusps). Let $S$ denote the set of singular points
(`corners') of $\d U$, and let $S_N$ be the small neighborhood of
$S$ given by
$$S_N=\left\{z\in \d U:
\dist (z,S)<\frac{b'\sqrt{\log N}}{\sqrtn}\right\} \;,$$ where
$b'>0$ is to be chosen below. We shall show that:
\begin{itemize} \item[i) ] \eqref{univ} holds uniformly for $z_0\in \d
U\sm S_N$;

\item[ii) ]
$ \displaystyle\sup_{z\in\d U} |\Upsilon^N_j(z)|= O\left(
N^{2j-m-1/2+\ep}\right),$  for $1\le j\le k$.\end{itemize}

Let us assume (i)--(ii) for now.  Since $\vol_{2m-1}S_N
=O\left(\frac{\sqrt{\log N}}{\sqrtn}\right)$, the estimate (ii)
implies that $$\int_{S_N} \Upsilon^N_j\wedge \om^{m-j} =
O(N^{2j-m-1+\ep})\;,$$ and hence by \eqref{VjU},
$$V_j^N(U)=\frac 1{(m-k)!^2}\left(\frac N\pi
\right)^{2k-2j} \int_{\d U\sm S_N}
\Upsilon^N_j\wedge\left[\om^{m-j} +O\left(\frac
1N\right)\right]+O\left(N^{2k-m-1+\ep}\right)\,.$$

It then follows from (i) and \eqref{Upj2} that \begin{eqnarray*}
V_j^N(U) &=& \frac {N^{2k-m-1/2}}{(m-k)!^2\,\pi^{2k-2j}}
\left[\int_{\d U\sm S_N} \Upsilon^\infty_j \wedge \om^{m-j}
+O(N^{-1/2+\ep})\right]\\ &=&\frac
{c_{mj}\,N^{2k-m-1/2}}{(m-k)!^2\,\pi^{2k-2j}}\left[\vol(\d U\sm
S_N)+O(N^{-1/2+\ep}) \right].\end{eqnarray*} Then by \eqref{VU}
\begin{eqnarray*} \var\big(\vol_{2m-2k}[Z_{s_1^N,\dots,s_k^N}\cap U]\big)
&=& N^{2k-m-1/2}\left[\nu_{mk}\,\vol_{2m-1}(\d U\sm S_N)
 +O(N^{-\frac 12 +\ep})\right]\\
&=& N^{2k-m-1/2}\left[\nu_{mk}\,\vol_{2m-1}(\d U)
 +O(N^{-\frac 12 +2\ep})\right]\;,\end{eqnarray*} which is our desired
formula.

It remains to prove (i)--(ii). To verify  (i), for each point
$z_0\in \d U\sm S$, we choose holomorphic coordinates $\{w_j\}$
and non-holomorphic coordinates $\{\wt w_j\}$ as above.  We can
choose these coordinates on a geodesic ball $V_{z_0}$ about $z_0$
of a fixed radius $R>0$ independent of the point $z_0$, but if
$z_0$ is near a corner, $\d U$ will coincide with $\{\Im \wt
w_1=0\}$ only in a small neighborhood of $z_0$. To be precise, we
let $D_{z_0}$ denote the connected component of $V_{z_0}\cap \d
U\sm S$ containing $z_0$. Then we choose $\phi
\in\ccal^2(V_{z_0})$ with $\phi(0)=0,\ d\phi(0)=0$, such that
\begin{equation}\label{Kz}\{w\in V_{z_0}:\Im w_1+\phi(w)=0\}
=\{\Im \wt w_1=0\}\supset D_{z_0}\;.\end{equation} We let
$$C= \sup _{z\in \d U\sm S}\; \frac{\dist (z,S)}{\dist (z,\d U\sm
D_z)}\ge 1\;.$$ Choose $N_0>0$ such that  $b\sqrt{\frac{\log
N_0}{N_0}} <R$; then
$$\textstyle\left\{w\in\d
U:\dist(z_0,w)<b\sqrt{\frac{\log N}{N}}\right\}\subset V_{z_0}\;,\quad
\mbox{for }\ N\ge N_0\;.$$

We recall that our assumption that $\d U$ is piecewise $\ccal^2$
{\it without cusps\/} means that $\overline U$ is locally
$\ccal^2$ diffeomorphic to a polyhedral cone, which implies that
$C<+\infty$.  We now let $b'=Cb$, where $b=\sqrt{2m+3}$ as before.

Consider any point $z_0\in \d U\sm S_N$, $N\ge N_0$.  Then
$$\dist (z_0,\d U\sm D_{z_0})\ge \frac {\dist (z_0,S)}C \ge
\frac{b'\sqrt{\log N}}{C\,\sqrtn} = \frac{b\sqrt{\log
N}}{\sqrtn}\;.$$ Thus by our far-off-diagonal decay estimate
\eqref{d2Qdecay}, the points in $\d U\sm D_{z_0}$ contribute
negligibly to the integral in \eqref{Upj}, so that integral can be
taken over the set $$\left\{|v|\le b\sqrt{\log N}:\,z_0+\frac
v\sqrtn\in D_{z_0}\right\}\,,$$ which is mapped by $\tau_N$ into
$\R\times \C^{m-1}$.  Then \eqref{univ0} holds, and \eqref{univ}
follows as before.

To verify  (ii),  we must show that the integral in the right side
of \eqref{Upj},
\begin{equation*}
\wt\Upsilon^N_j(z_0):=\int_{\left\{|v|\le b\sqrt{\log
N}:\,z_0+\frac v\sqrtn\in\d U\right\}}\textstyle\frac 14\,
F''(\half|v|^2)\,   \dbar(\bar z\cdot v)
\wedge \dbar|v|^2  \wedge
\left(\Var^{z_0}_\infty\right)^{j-1}\wedge\left(\frac i2\ddbar
|v|^2\right)^{m-j}\;,\end{equation*} is $O(N^\ep)$ uniformly for
$z_0\in \d U$. By Lemma \ref{d4Qas}, $\Var^{z_0}_\infty( v) =
O(|v|^{-2})$. Furthermore, $$\left| \frac 1{e^{| v|^2}-1}\, (
v\cdot d\bar z)\wedge ( v\cdot d\bar{ v})\right| \le \frac{\sqrt
m\, |v|^2} {e^{|v|^2} -1}\le \sqrt m$$ (using Euclidean norms in
the $z$ and $v$ variables ), and hence
\begin{equation}\label{showii} |\wt\Upsilon^N_j(z_0)| \le A_{jm}
\int_{\left\{|v|\le b\sqrt{\log N}:\,z_0+\frac v\sqrtn\in\d
U\right\}} |v|^{-2j+2} \,d\vol_{2m-1}^E(v)\;,\end{equation} for
universal constants $A_{jm}$,where $\vol^E$ denotes Euclidean
volume.  Rewriting \eqref{showii} in terms of the original
variables $w=z_0+\frac v\sqrtn$, we have
\begin{equation*} |\wt\Upsilon^N_j(z_0)| \le A_{jm}\,N^{-j+m +1/2}
\int_{\left\{w\in\d U:\,|w-z_0|\le b\sqrt{\frac{\log N}N}\right\}}
|w-z_0|^{-2j+2} \,d\vol_{2m-1}^E(w)\;.\end{equation*}

For each point $P\in\d U$, we choose a closed neighborhood $V_P\in
M$ of $P$ and a $\ccal^2$ diffeomorphism $\rho_P:V_P\to \R^{2m}$
mapping $V_P\cap\d U$ to the boundary of a polyhedral cone
$K_P\subset \R^{2m}$.  Then for $N$ sufficiently large, for all
$z_0\in \d U$, the set $\left\{w\in\d U:\,|w-z_0|\le
b\sqrt{\frac{\log N}N}\right\}$ is contained in one of the  $V_P$.
We then make the (nonholomorphic) coordinate change $\wt w=
\rho_P(w)$.  Since the diffeomorphisms $\rho_P$ have bounded
distortion, we then have
\begin{equation}\label{iia} |\wt\Upsilon^N_j(z_0)| \le A'_{jm}\,N^{-j+m +1/2}
\int_{\left\{\wt w\in\d K_P:\,|\wt w-\wt z_0|\le b\sqrt{\frac{\log
N}N}\right\}} |\wt w-\wt z_0|^{-2j+2} \,d\vol_{2m-1}^E(\wt
w)\;.\end{equation}  Let $q\in\Z^+$ be the maximum number of
facets in $\d K_P$.  We easily see that
\begin{eqnarray}\label{iib} \int_{\left\{\wt w\in\d K_P:\,|\wt
w-\wt z_0|\le b\sqrt{\frac{\log N}N}\right\}} |\wt w-\wt
z_0|^{-2j+2} \,d\vol_{2m-1}^E(\wt w) & \le & q
\int_{\left\{x\in\R^{2m-1}:\,|x|\le b\sqrt{\frac{\log
N}N}\right\}} |x|^{-2j+2} \,dx\nonumber\\ &=& \mbox{const.}\times
\left(\frac{\log N}N\right)^{m-j-1/2}\;.
\end{eqnarray}  Combining \eqref{iia}--\eqref{iib}, we conclude
that $\wt\Upsilon^N_j(z_0)=O(N^\ep)$   and thus (ii) holds, which verifies
the formula of Theorem \ref{volume} for the general case
where $\d U$ has corners.

\subsection{Positivity of the constants $\nu_{mk}$}\label{positivity}

To complete the proof of Theorems \ref{number} and \ref{volume}, we must show
that the leading coefficients
$\nu_{mk}$  are positive.  In fact, we shall show that each
of the coefficients $c_{mj}$ (defined by \eqref{Upj2}) is positive, and then
the positivity of the $\nu_{mk}$ follows from \eqref{numk}.

We begin  by computing  the coefficient $\nu_{m1}$ in the
codimension-one case: By
\eqref{Upj1},
\begin{eqnarray}\Upsilon^\infty_1(z_0)&=& \frac {1}{4\pi^2}
\sum_{j,k}\left[\int_{\R\times \C^{m-1}} \frac{v_jv_k d\bar v_k}
{e^{|v|^2}-1} \wedge \left(\textstyle\frac i2\ddbar
|v|^2\right)^{m-1}\right] d\bar z_j\nonumber
\\&=& \frac {(m-1)!}{4\pi^2 }\sum_{j=1}^m
\left[\int_{\R\times \C^{m-1}}\frac{v_j v_1} {e^{|v|^2}
-1}\,d\vol_{\R\times \C^{m-1}}(v)\right] d\bar
z_j\,.\label{Up2}\end{eqnarray}

By \eqref{VU}--\eqref{VjU} with $k=1$ and \eqref{univ}, we have
\begin{equation}\label{V1U}
\var\big(\vol_{2m-2k}[Z_{s_1^N}\cap U]\big) = \frac
{N^{3/2-m}}{(m-1)!^2}\left[ \int_{\d U}
\Upsilon^\infty_1(z)\wedge\om(z)^{m-1}
+O\left(N^{-1/2+\ep}\right)\right]\,.\end{equation} Since
$d\vol_{\d U}(z_0)= dx_1\wedge \frac 1{(m-1)!}\om(z_0)^{m-1}$,
only the $j=1$ term in \eqref{Up2} contributes to the integral in
\eqref{V1U}, and we  then have
\begin{equation}\var\big(\vol_{2m-2k}[Z_{s_1^N}\cap
U]\big) = N^{3/2-m}\left[ \int_{\d U} \nu_{m1}d\vol_{\d U}
+O\left(N^{-1/2+\ep}\right)\right],\end{equation} where
\begin{eqnarray*} \nu_{m1} &=& \frac 1{4\pi^2}
\int_{\R^{2m-1}} \frac{v_1^2}{e^{|v|^2}-1}\,dv\ =\ \frac 1
{4\pi^2(2m-1)} \int_{\R^{2m-1}} \frac{|v|^2}{e^{|v|^2}-1}\,dv
\\&=&  \frac 1 {4\pi^2(2m-1)} \,
\frac{2\pi^{m-1/2}}{\Gamma(m-1/2)}\int_0^\infty
\frac{r^{2m}}{e^{r^2}-1}\,dr
\\&=& \frac {\pi^{m-5/2}}{4\,\Gamma(m+1/2)}\sum_{k=1}^\infty
\int_0^\infty e^{-kr^2}\, r^{2m}\,dr
\\&=&  \frac {\pi^{m-5/2}}{4\,\Gamma(m+1/2)}\sum_{k=1}^\infty
\frac {\Gamma(m+1/2)}{2\,k^{m+1/2}} \ =\
\frac{\pi^{m-5/2}}{8}\;\zeta\Big(m+\half\Big)\;,
\end{eqnarray*} as stated in the theorem.

  We now determine the $c_{mj}$, for $1\le j\le m$. By
\eqref{Upj1s}--\eqref{Upj2}, 
\begin{multline}\label{p1} 
c_{mj}
\;dx_1\wedge \textstyle (\frac i2 dz_2\wedge d\bar
z_2)\wedge\cdots\wedge (\frac i2 dz_m\wedge d\bar z_m)\ =\ c_{mj}\,d\vol_{\d
U}\\ =\ \frac
{(-1)^{j-1}}{(4\pi^2)^j}\int_{\R\times
\C^{m-1}}\frac 1{(e^{|v|^2}-1)^j}\, (v\cdot d\bar z)\wedge (v\cdot d\bar v)
\hspace{2.2in}
\\ \wedge\big[(d\bar z\cdot dv)\wedge(dz\cdot d\bar v)+(dz\cdot d\bar z)\wedge
((dv\cdot d\bar v)\big]^{j-1} \textstyle \wedge \big[\frac i2
\ddbar|z|^2 \wedge \frac i2 \ddbar|v|^2\big]^{m-j},
\end{multline}
where $v\cdot d\bar z = \sum v_\al\,d\bar z_\al,\ d\bar z\cdot dv
=\sum d\bar z_\al\wedge dv_\al$, etc.  As in the codimension one case, only
the term $v_1^2 d\bar z_1\wedge d\bar v_1$ in $ (v\cdot d\bar z)\wedge (v\cdot
d\bar v)$ contributes to the integral, and we obtain
\begin{multline}\label{p2} 
c_{mj}
\,d\vol_{\d U}
=  \left(\frac i2\right)^{2m-2}\frac
{1}{4\pi^{2j}}\int_{\R\times
\C^{m-1}}\frac {v_1^2}{(e^{|v|^2}-1)^j}\, d\bar z_1\wedge d\bar v_1
\wedge(A+B)^{j-1}\wedge B^{m-j}\,,
\end{multline}
where $$\begin{array}{llll}  A &=&\di  \sum_{\al,\be =2}^m A_{\al\be}\,,
\quad & A_{\al\be} = dz_\al\wedge d\bar z_\be \wedge dv_\be\wedge d\bar
v_\al\,,\\
B &=&\di  \sum_{\al,\be =2}^m B_{\al\be}\,,
\quad & B_{\al\be} = dz_\al\wedge d\bar z_\al \wedge dv_\be\wedge d\bar
v_\be\,.\end{array}$$

Writing \begin{equation}\label{p3}(A+B)^{j-1}\wedge B^{m-j} =
\sum_{\ell=0}^{j-1} {j-1\choose \ell}A^\ell\wedge B^{m-\ell-1}\,,\end{equation}
we see that it suffices to show that $\left(\frac i2\right)^{2m-2}\,
A^\ell\wedge B^{m-\ell-1}$ is positive for $0\le \ell \le m-1$.  We have
\begin{equation}\label{product} A^\ell\wedge B^{m-\ell-1} = \sum
_{\al_1,\be_1=2}^m \cdots \sum
_{\al_{m-1},\be_{m-1}=2}^m A_{\al_1\be_1}\wedge \cdots\wedge 
A_{\al_\ell\be_\ell}\wedge  B_{\al_{\ell+1}\be_{\ell+1}}\wedge \cdots\wedge 
B_{\al_{m-1}\be_{m-1}}\;.\end{equation}
We claim that each term in \eqref{product}  equals  $B_{22}\wedge
B_{33}\wedge\cdots\wedge B_{mm}$ if $\al_1,\dots,\al_{m-1}$ is a 
permutation of $ 2,\dots,m$, and 
$\be_1,\dots,\be_{m-1}$ is obtained from $\al_1,\dots,\al_{m-1}$ by permuting
$\al_1,\dots,\al_\ell$ and also permuting $\al_{\ell+1},\dots,\al_{m-1}$. 
Otherwise, the term clearly  vanishes. To verify the claim, by simultaneously
permuting the coordinates
$\{z_\al\}$ and $\{v_\al\}$, we can assume without loss of generality that
$\al_p=p+1$ ($1\le p\le m-1$); starting with
$\be_p=\al_p$, where the claim is a tautology, we  note that a transposition in
$\be_1,\dots,
\be_\ell$
 transposes  two $(d\bar
v_\be \wedge dz_\be)$'s, while a transposition in
$\be_{\ell+1},\dots,\be_{m-1}$ transposes  two $(d
v_\be \wedge d\bar v_\be)$'s, leaving the sign unchanged.

If follows that 
\begin{equation}\label{p5} \textstyle \left.\left(\frac i2\right)^{2m-2}
d\bar z_1\wedge d\bar v_1 \wedge  A^\ell\wedge
B^{m-\ell-1}\right|_{\d U\times (\R\times \C^{m-1})} =
(m-1)!\,\ell!\,(m-\ell-1)!
\,d\vol_{\d U}
\wedge d\vol_{\R\times
\C^{m-1}},\end{equation} and thus  the  $c_{mj}$ are positive,
completing the proof of Theorem \ref{volume}. \qed

\bigskip Combining \eqref{numk}, \eqref{p2}, \eqref{p3} and \eqref{p5}, we
obtain the explicit formula:

\begin{equation}\label{nu}
\nu_{mk}= \frac{\pi^{m-2k-1/2}\,k!\,(m-1)!}{4\,\Gamma(m+1/2)\,(m-k)!^2}\, 
\sum_{j=1}^k \frac{n(m,j)}{j\,(k-j)!}\,
\int_0^\infty \frac{r^{2m}\,dr}{(e^{r^2}-1)^j}\,,\end{equation}
where $n(m,j)\in\Z^+$ is given by \begin{equation}\label{nu1}n(m,j) =
\sum_{\ell=1}^{j}
\frac{(m-\ell)!}{(j-\ell)!}\;.\end{equation}

\section{ Appendix: Proof of Theorem \ref{near-far}}\label{s-proof}

In this
appendix, we sketch the proof of the off-diagonal \szego
asymptotics theorem. The argument is essentially contained in
\cite{SZ2}, but we add some details relevant to the estimates in
Theorem \ref{near-far}.

The \szego kernels $\Pi_N(x,y)$ are the Fourier coefficients of
the total
\szego projector $\Pi(x,y):\lcal^2(X)\to \hcal^2(X)$; i.e.
$\Pi_N(x,y)=\frac 1{2\pi}\int
e^{-iN\theta}\Pi(e^{i\theta}x,y)\,d\theta$. The estimates for
$\Pi_N(z,w)$ are then  based on the  Boutet de
Monvel-Sj\"ostrand construction of an  oscillatory integral
parametrix for the \szego kernel:
\begin{equation}\label{oscint}\begin{array}{c}\Pi (x,y) =
S(x,y)+E(x,y)\;,\\[12pt] \mbox{with}\;\; S(x,y)=
\int_0^{\infty} e^{i t \psi(x,y)} s(x,y,t ) dt\,, \qquad
E(x,y)\in
\ccal^\infty(X \times X)\,.\end{array}
\end{equation} The amplitude has the form $s \sim \sum_{k =
0}^{\infty} t^{m -k} s_k(x,y)\in S^m(X\times X\times \R^+)$. The
phase function $\psi$ is of positive type, and as described in
\cite{BSZ2},  is given by:
\begin{equation}
\psi(z,
\theta, w, \phi) = {i} \left[1 - \frac{a(z,\bar{w})}{
\sqrt{a(z)} \sqrt{a(w)}}\; e^{i (\theta - \phi)}\right]\;,
\label{psiphase}\end{equation} where
$a\in\ccal^\infty(M\times M)$ is an almost holomorphic
extension  of the function
$a(z,\bar z):=a(z)$ on the anti-diagonal $A=\{(z,\bar
z):z\in M\}$, i.e.,\ $\dbar a$ vanishes to infinite order
along $A$. We recall from
\eqref{a} that $a(z)$ describes the Hermitian metric on $L$ in
our preferred holomorphic frame at
$z_0$, so by \eqref{preferred}, we have $
a(u)=1+|u|^2+O(|u|^3)$, and hence
\begin{equation}\label{auv} a(u,\bar v)= 1+ {u\cdot\bar v}
+O(|u|^3+|v|^3)\;.\end{equation} For further background and
notation on complex Fourier integral operators we refer to
\cite{BSZ2} and to the original paper of Boutet de Monvel and
Sj\"ostrand
\cite{BS}.

As above,  denote the $N$-th Fourier coefficient of these
operators relative to the $S^1$ action by $\Pi_N = S_N + E_N$.
Since $E$ is smooth, we have $E_N(x,y) = O(N^{-\infty})$, where
$O(N^{-\infty})$ denotes a quantity which is  uniformly
$O(N^{-k})$ on $X\times X$ for all positive $k$. Then,
$E_N(z,w)$ trivially satisfies the
remainder estimates in Theorem \ref{near-far}.

Hence it is only necessary to verify that the oscillatory integral
\begin{equation}S_N(x,y)  =   \int_0^{2\pi} e^{- i
N \theta}  S( e^{i\theta} x,y) d\theta  =  \int_0^{\infty}
\int_0^{2\pi}  e^{- i N \theta+it  \psi( e^{i\theta} x,y)}
s(e^{i\theta} x,y,t) d\theta dt \end{equation} satisfies Theorem
\ref{near-far}. This follows from an analysis of the  stationary
phase method  and remainder estimate for the rescaled parametrix
 \begin{equation}\begin{array}{l}\displaystyle S^{z_0}_N\left(
\frac{u}{\sqrt{N}}, 0;  \frac{v}{\sqrt{N}}, 0\right)
= N \int_0^{\infty} \int_0^{2\pi}
 e^{ i N \left( -\theta + t\psi\big(  \frac{u}{\sqrt{N}}, \theta;
\frac{v}{\sqrt{N}}, 0\big)\right)} s\left( \frac{u}{\sqrt{N}},
\theta; \frac{v}{\sqrt{N}}, 0, Nt\right) d\theta dt
\,,\end{array}\label{SN}\end{equation} where we changed  variables
$t \mapsto N t$. For background on the stationary phase method
when the phase is complex we refer to \cite{H}. We are
particularly interested in the dependence of the stationary phase
expansion and remainder estimate on the parameters $(u, v)$
satisfying the constraints in (i)-(ii) of Theorem \ref{near-far}.

To clarify the constraints, we recall  from  \cite[(95)]{SZ2}  that
the \szego kernel satisfies the following far from diagonal
estimates:
\begin{equation} \label{CONSTRAINTS}  
\left|\nabla^j_h \Pi_N(z,w)\right|=O(N^{-K})\qquad \mbox{for
all } \;j,  K \; \mbox{when} \;\;
\dist(z,w)\geq\,\frac{N^{1/6}}{\sqrt{N}} \;. \end{equation}  Hence we
may assume from now on that $z = z_0 + \frac{u}{\sqrt{N}}, w = z_0
+ \frac{v}{\sqrt{N}}$ with
\begin{equation} \label{CONSTRAINTS2} |u|+|v|\le \delta N^{1/6}
\end{equation}
for a sufficiently small constant $\delta > 0$.

By \eqref{psiphase}--\eqref{auv},  the rescaled phase in \eqref{SN} has the
form:
\begin{equation}
 \wt\Psi:= t \psi\left( \frac{u}{\sqrt{N}},\theta;  \frac{
v}{\sqrt{N}}, 0\right)  -\theta  = it \left[ 1 - \frac{a\left(
\frac{u}{\sqrt{N}},  \frac{\bar v}{\sqrt{N}}\right)}{a\left(
\frac{u}{\sqrt{N}}, \frac{\bar u}{\sqrt{N}}\right)^{\frac 12}
a\left( \frac{v}{\sqrt{N}},  \frac{\bar v}{\sqrt{N}}\right)^{\frac
12}}\;  e^{i \theta}\right] -  \theta
\end{equation}
and the  $N$-expansion
\begin{equation}\label{entirephase} \wt\Psi=it[ 1 - e^{i \theta}]
-  \theta -\frac{it}{N}\psi_2(u,v) e^{i \theta} + t
R_3^\psi(\frac{u}{\sqrtn},\frac{v}{\sqrtn}) e^{i \theta}\,,
\end{equation}
where
\begin{equation*} \psi_2(u,v) =
u \cdot\bar{v} - \half(|u|^2 + |v|^2)=-\half|u-v|^2+i\,\Im(u\cdot\bar
v)\;.\end{equation*} After
multiplying by $iN$, we move the last two terms of \eqref{entirephase}
into the amplitude.
 Indeed, we    absorb all of  $ \exp\{(\psi_2 +i
NR_3^\psi)te^{i \theta}\}$   into the amplitude so that (\ref{SN})
is an oscillatory integral
\begin{equation}\label{phase-amplitude}  N
\int_0^{\infty} \int_0^{2\pi}
e^{iN\Psi(t,\theta)}A(t,\theta;z_0,u,v)d\theta dt + O(N^{-\infty})
\end{equation}with phase
\begin{equation}\label{phase}\Psi(t,\theta): = it ( 1 -
e^{i\theta})- \theta\end{equation} and with amplitude
\begin{equation}\label{amplitude}
A(t,\theta;z_0,u,v):=  e^{ t e^{i \theta} \psi_2(u,v)  + it e^{i
\theta} N R_3^\psi( \frac{u}{\sqrt{N}},\frac{ v}{\sqrt{N}})}\,
s\big( \frac{u}{\sqrt{N}}, \theta; \frac{v}{\sqrt{N}}, 0, Nt\big).
\end{equation}

 The  phase $\Psi$ is independent of the
parameters $(u,v)$, satisfies $\Re (i \Psi) = - t (1 - \cos
\theta) \leq 0$ and has a unique critical point at $ \{t=1, \theta
= 0\}$ where it vanishes.

The factor $e^{ t e^{i \theta} \psi_2(u,v) }$ is of exponential
growth
 in some regions.   However, since it is  a
 rescaling of a complex phase of positive type, the complex phase
 $i N \Psi +  t e^{i \theta} \psi_2(u,v)$ is of positive type,
 \begin{equation} \label{NEG} \Re (i N \Psi +  t e^{i \theta} \psi_2(u,v))
   <
 0\end{equation}
 once the cubic remainder $N t e^{i \theta}
R_3^\psi(\frac{u}{\sqrtn},\frac{v}{\sqrtn})$ is smaller than $i N
\Psi + t e^{i \theta} \psi_2(u,v)$, which occurs
 for all $(t, \theta, u, v)$ when $(u, v)$ satisfy
 (\ref{CONSTRAINTS2}) with $\delta$ sufficiently small.

 To estimate the
 joint rate of decay in $(N, u, v)$, we  follow the  stationary
phase expansion and remainder estimate in Theorem 7.7.5 of
\cite{H}, with extra attention to    the unbounded  parameter $u$.

The first step is to use a  smooth partition of unity $\{\rho_1(t,
\theta),\rho_2(t, \theta )\}$ to decompose the integral (\ref{SN})
into a region  $(1 - \epsilon, 1 + \epsilon)_t \times (-\epsilon,
\epsilon)_{\theta}$  containing the critical point and one over
the complementary set  containing no critical point. We claim that
the  $\rho_2$ integral  is of order $N^{-\infty}$ and can be
neglected. This follows by repeated partial integration as in the
standard proof together with the fact that the exponential factors
in (\ref{NEG}) decay, so that the estimates are integrable and
uniform in $u$.

We then apply  \cite{H}  Theorem 7.7.5 to  the $\rho_1$ integral.
 The first term
of the stationary phase expansion equals $N^m e^{ t e^{i \theta}
\psi_2(u,v)}$ and the remainder satisfies

\begin{equation} \label{remainder} |\wh R_{J}(P_0, u, v, N)| \leq C
N^{- m + J} \sum_{|\alpha| \leq 2J+2} \sup_{t, \theta}
|D^{\alpha}_{t, \theta} \rho_1 A(t,\theta;P_0,u,v) |.
\end{equation}
From the formula in (\ref{amplitude}) and the  fact that $s$ is a
symbol, $A$ has a polyhomogeneous expansion of the form
\begin{eqnarray} A(t,\theta;P_0,u,v) &=& \rho_1(t, \theta)   
e^{ t e^{i \theta} \psi_2(u,v)}
 N^m \left[\sum_{n = 0}^{K} N^{-n/2} f_{n}(u, v;t, \theta,
P_0)+R_K(u, v, t,\theta)\right]\,,\nonumber\\
&&\qquad\qquad |\nabla^jR_{Nk}(u,v)|\le C_{jk\ep b}e^{\ep (|u|^2 +
|v|^2)} N^{-\frac{K+1}{2}}\big). \label{Alarge}\end{eqnarray} The
exponential remainder factor $e^{\ep (|u|^2 + |v|^2)}$  comes from
the fact   $\Re e^{i \theta} \psi_2 = \cos \theta \Re \psi - \sin
\theta \Im \psi$ with $\Re \psi \leq 0$ and  $|\sin \theta| <
\epsilon$ on the support of $\rho_1$. Hence, the supremum of the
amplitude in a neighborhood of the stationary phase set (in the
support of $\rho_1$) is bounded by $e^{\epsilon |\Im \psi_2|}$.
The remainder term is smaller than the main term asymptotically as
$N \to \infty$ as long as $(u,v)$ satisfies (\ref{CONSTRAINTS2}).
Part(i) of Theorem \ref{near-far} is an immediate consequence of
(\ref{Alarge}) since
$e^{\ep (|u|^2 + |v|^2)} \leq N^{\epsilon}$ for $|u| + |v|
\leq \sqrt{\log N}$.

To prove part (ii), we may assume from
(\ref{CONSTRAINTS})--(\ref{CONSTRAINTS2}) that $\sqrt{\log N} \leq
|u|+|v|  \leq\de\,N^{1/6}$. In this range the
asymptotics (\ref{Alarge}) are valid. We first rewrite the horizontal
$z$-derivatives $\frac{\d^h}{\d z_j}$ as $u_j$ derivatives, which for
$L^N$ have the form  $\sqrt{N}
\frac{\d}{\d u_j} - N A_j(\frac u{\sqrtn})$ and thus $\nabla_h$
contributes a factor of
 $\sqrtn$. We thus obtain an
asymptotic expansion and remainder for $\nabla^j_h \Pi_N(z,w)$ by
applying $\nabla^j_h$ to the expansion (i) with $k = 0$: $$
\Pi^\Heis_1(u,\theta;v,\phi)\left[1+ N^{-1/2} R_{N0}(u,v)\right].
$$  The operator $\nabla^j_h $ contributes a factor of $N^{j/2}$
to each term, and thus
\begin{eqnarray*}\left|\nabla^j_h
\Pi_N(z,w)\right|& = & O\left(N^{m + j/2} \, e^{-(1 - \epsilon)
\frac{|u|^2 + |v|^2}{2}}\right)  \\
& = & O(N^{ -k })\qquad \mbox{uniformly for }\;\;
{|u|^2 + |v|^2} \ge (j+2k+2m+\ep'){ \log N} \;,
\end{eqnarray*} where $\ep'=(j+2k+2m+1)\ep$.
\qed

\end{document}